\documentclass[oneside]{article}
\usepackage[intlimits,reqno]{amsmath}
\usepackage{amsthm, graphics}
\usepackage{amssymb}
\usepackage{verbatim}
\include{diagxy}
\input{xy}
\xyoption{all}

\usepackage{enumerate}
\usepackage{amsmath}
\usepackage{verbatim}
\usepackage{bbm}
\newcommand{\norm}[1]{\left\Vert#1\right\Vert}
\renewcommand{\sc}[2]{\langle #1|#2 \rangle}

\newcommand{\proj}[2]{|#1 \rangle \langle #2|}

\newcommand{\bra}[1]{\left\langle #1 \right|}
\newcommand{\ket}[1]{\left| #1 \right\rangle}

\newcommand{\A}{\mathcal{A}}

\newcommand{\G}{\mathcal{G}}

\newcommand{\U}{\mathcal{U}}

\renewcommand{\H}{\mathcal{H}}

\renewcommand{\ll}{l}
\newcommand{\rr}{r}
\newcommand{\Ll}{\mathcal{L}}
\newcommand{\Ss}{\textbf{s}}
\newcommand{\Tt}{\textbf{t}}

\newcommand{\tto}{\rightrightarrows}

\newcommand{\M}{\mathfrak{M}}

\renewcommand{\to}{\rightarrow}
\newcommand{\prf}[1]{\begin{proof}#1\end{proof}}

\newcommand{\ol}{\overline}

\newcommand{\be}{\begin{equation}}
\newcommand{\ee}{\end{equation}}
\newcommand{\bse}{\begin{subequations}}
\newcommand{\ese}{\end{subequations}}
\newcommand{\ben}{\begin{enumerate}}
\newcommand{\een}{\end{enumerate}}
\newcommand{\bit}{\begin{itemize}}
\newcommand{\eit}{\end{itemize}}
\newcommand{\bex}{\begin{example}}
\newcommand{\eex}{\begin{flushright}$\diamondsuit$\end{flushright}\end{example}}

\newtheorem{thm}{Theorem}[section]
\newtheorem{cor}[thm]{Corrolary}

\newtheorem{prop}[thm]{Proposition}

\theoremstyle{definition}

\newtheorem{example}{Example}
\numberwithin{example}{section}
\theoremstyle{remark}

\numberwithin{equation}{section}

\DeclareMathOperator{\pr}{pr}

\renewcommand{\emph}[1]{{\bf #1}}

\begin{document}\large{
\title{Algebroids associated to the groupoid of partially invertible elements of a $W^*$-algebra}
\author{Anatol Odzijewicz,\ Grzegorz Jakimowicz,\ Aneta Sli\.{z}ewska \\Institute of Mathematics\\
University in Bia{\l}ystok
\\Akademicka 2, 15-267 Bia{\l}ystok, Poland}

\maketitle
\begin{abstract}
In the paper we study the algebroid $\A(\M)$ of the groupoid $\G(\M)\tto\Ll(\M)$ of partially invertible elements over the lattice $\Ll(\M)$  of orthogonal projections of a $W^*$-algebra $\M$. In particular the complex analytic manifold structure of these objects is investigated.
The expressions on the algebroid Lie brackets for $\A(\M)$ and related algebroids are given in noncommutative
operator coordinates in the explicit way. We also prove statements describing structure of the groupoid of partial isometries $\U(\M)\tto\Ll(\M)$ and the frame groupoid $\G^{lin}\A_{p_0}(\M)\tto\Ll_{p_0}(\M)$ as well as their algebroids.
\end{abstract}

\tableofcontents
\vspace{1cm}

\section{Introduction}
Nowadays the theory of Lie groupoids and Lie algebroids is an important and invaluable part of contemporary differential geometry. One can consider Lie groupoids and Lie algebroids as a framework for the investigation of symmetry of objects having fibre bundle structure. They  also are a natural generalizations of Lie groups and Lie algebras, respectively. In the monograph \cite{mac} of Mackenzie, whose personal contribution to the subject is notable, one finds compact presentation of the subject as well as historical references at the end of each chapter.

While the concept of groupoid in topology \cite{brown} and Lie groupoids and Lie algebroids in differential geometry \cite{pradines,pradines1,pradines2} appeared in 60s and 70s of the last century the growth of interest to these notions was truly inspired by problems of mathematical physics. It was firstly Poisson geometry, where after the paper of Karasev \cite{karasev} and seminar note of Coste, Dazord  and Weinstein \cite{coste} the symplectic realization of Poisson manifold by symplectic groupoid was defined and investigated. Then it was observed that the existence of Poisson structure on a manifold $M$ is equivalent to Lie algebroid structure on its cotangent bundle $T^*M$. Subsequently many other researchers, first of all Weinstein  with his collaborators used groupoid and algebroid methods in mechanics and in problems of quantization, see Cannas da Silva, Weinstein \cite{Wei}, Connes \cite{conne} and references of therein.

Admitting some modifications of the basic definitions one can investigate the above mentioned structures in the infinite dimentional case, i.e. in the framework of the category of smooth Banach manifolds. An example of the investigations of Banach Lie algebroids one finds in \cite{pellet}. However, a crucial difficulty appears when the modeling Banach space does not have Schauder basis, e.g.  it happens for Banach algebra of bounded operators $L^\infty(\H)$ of Hilbert space $\H$. The reason is that in this case does not exist coordinate description of the investigated structures, in particular case the Poisson structure.

Nevertheless the category of $W^*$-algebras (von Neumann algebras) is the one which generates a rich class of Banach-Lie groupoids and Banach-Lie algebroids as well as the related class of Banach Poisson manifolds which have nice properties and could be handled in the operator coordinate manner which is strictly associated to $W^*$-algebra structure. We were motivated in our investigations by importance of von Neumann algebras in quantum physics and the fact that the predual $\M_*$ of the $W^*$-algebra $\M$ has canonically defined structure of Banach Lie Poisson space (see \cite{OR}) and moreover the precotangent bundle $T_*G(\M)$ of the Banach Lie group $G(\M)$ of invertible elements of $\M$ is a weak symplectic realization of $\M_*$.

The main place in our considerations is occupied by the Banach Lie groupoid $\G(\M)\tto \Ll(\M)$ of partially invertible elements of $\M$ with the lattice of orthogonal projections $\Ll(\M)$ as the base manifold which was defined and described in \cite{OS}.

In Section 2 we investigate the family of locally trivial transitive subgroupoids $\G_{p_0}(\M)\tto \Ll_{p_0}(\M)$ of the groupoid $\G(\M)\tto \Ll(\M)$, parametrized by $p_0\in \Ll(\M)$. Namely, we define a complex analytic atlas on $\G_{p_0}(\M)\tto \Ll_{p_0}(\M)$ consistent with its groupoid structure, i.e. we show that $\G_{p_0}(\M)\tto \Ll_{p_0}(\M)$ is a complex analytic groupoid.

In Section 3 we describe the Banach-Lie algebroid  $\A_{p_0}(\M)$ of the groupoid $\G_{p_0}(\M)\tto \Ll_{p_0}(\M)$. Our description  is based on the groupoid isomorphism given in Proposition \ref{prop22}, which shows that one can consider the groupoid $\G_{p_0}(\M)\tto \Ll_{p_0}(\M)$ as the gauge groupoid of the $G_0$-principal bundle with $P_0:=\G_{p_0}(\M)\cap \M p_0$ as the total space and the group $G(p_0\M p_0)$ of invertible elements of $p_0\M p_0$ denoted by $G_0$. Among others we present coordinate expressions (\ref{bracket}) and (\ref{nawias}) for Lie algebroid bracket of sections of $\A_{p_0}(\M)$. In Section 3 we also present detailed description of Banach-Lie groupoid $\U_{p_0}(\M)\tto \Ll_{p_0}(\M)$ of partial isometries as well as the corresponding  algebroid $\A^{u}_{p_0}(\M)$, see Proposition \ref{prop1U} and Proposition \ref{nawias4}.

As an example  we consider in Section 4 the algebroid of the frame groupoid of tautological vector bundle $\mathbb{E}\to G(N,\H)$ over Grassmannian of $N$-dimensional subspace of  the Hilbert space $\H$.

In Section 5 we discuss various groupoids associated in a canonical way to $\G_{p_0}(\M)\tto \Ll_{p_0}(\M)$ including the frame groupoid $\G^{lin}\A_{p_0}(\M)\tto \Ll_{p_0}(\M)$ of the algebroid $\A_{p_0}(\M)$.

Finally Section 6 contains the description of the algebroid $\A^{lin}_{p_0}(\M)$ which is the algebroid of the groupoid
$\G^{lin}\A_{p_0}(\M)\tto \Ll_{p_0}(\M)$ proving that $\A^{lin}_{p_0}(\M)$ is a subalgebroid of the algebroid  $\mathfrak{D}(TP_0)$ of derivations of  $\Gamma^\infty (TP_0)$. 

At the end the following two facts are worth to be noted. In the case $\M=L^\infty(\H)$ one can consider $\G_{p_0}(\M)\tto \Ll_{p_0}(\M)$, $\ \U_{p_0}(\M)\tto \Ll_{p_0}(\M)$, $\ \A_{p_0}(\M)$, $\ \A^u_{p_0}(\M)$, $P_0$ and $P_0^u$ as an universal objects in the corresponding categories. However, we will not discuss this question, leaving  it for subsequent paper. The groupoids $\G_{p_0}(\M)\tto \Ll_{p_0}(\M)$ as well as  the corresponding algebroids $\A_{p_0}(\M)$ provide interesting examples of the complex analytic Banach-Lie groupoids and algebroids, respectively.

\section{Groupoid of partially invertible elements of $W^*$-algebra}
Such a class of Banach-Lie groupoids was introduced and investigated in \cite{OS}.
Here we recall some necessary notions and statements concerning the
subject. By definition the groupoid $\G(\M)$ of partially invertible elements of
$W^*$-algebra $\M$ will consist of such elements $x\in\M$ for which
 $|x|=(x^*x)^{\frac{1}{2}}$ is an invertible element of the $W^*$-subalgebra
$p\M p\subset \M$, where $p$ is the support of $|x|$. We have natural
maps $\ll:\G(\M)\to \Ll(\M)$ and $\rr:\G(\M)\to \Ll(\M)$  of $\G(\M)$ on the complete lattice
$\Ll(\M)$ of orthogonal projections of the $W^*$-algebra $\M$, defined as the left and right support of $x\in \G(\M)$, respectively. Taking
$\Ll(\M)$ as the base set of $\G(\M)$ one can identify $\ll$ with the target map,
$\rr$ with the source map and inclusion $\varepsilon:\Ll(\M)\to \G(\M)$ with the
object inclusion map. The partial multiplication of $x,y\in \G(\M)$ is the
algebraic product in $\M$. Note that $\G(\M)\subset\M$ and $xy\in \G(\M)$ if
$\rr(x)=\ll(y)$. The two sided inversion $x^{-1}$ of $x\in\G(\M)$ is defined by its
polar decomposition \be\label{polar} x=u|x|\ee in the following way \be x^{-1}=|x|^{-1}u^*.\ee

One easily verifies that the above maps and operations define groupoid
structure $\G(\M)\tto \Ll(\M)$ on the set  $\G(\M)$  over the lattice $\Ll(\M)$.
\begin{prop} An element $x\in\M$ of $W^*$-algebra $\M$ belongs to $\G(\M)$ if and only if there exists $y\in\M$ such that
\be\label{y} yx=\rr(x)\qquad { and}\qquad xy=\ll(x)\ee
and one has
$\ll(y)=\rr(x)$ and $\rr(y)=\ll(x).$ The element $y$ belongs to $\G(\M)$ and is defined uniquely  by $x\in\G(\M)$.\end{prop}
\prf{ If $x\in\G(\M)$ then  $y=|x|^{-1}u^*$ satisfies relations (\ref{y}). Now, let $y\in\M$ satisfy conditions (\ref{y}). Then
\be\label{y1} u^*u=\rr(x)=yx=yu|x|\ee
and
\be\label{y2} uu^*=\ll(x)=xy=u|x|y.\ee From (\ref{y1}) we obtain \be\label{y3} \rr(x)=|x|yu.\ee From  (\ref{y1}) and (\ref{y3}) we see that $yu\in \rr(x)\M \rr(x)$ and $yu=|x|^{-1}$, i.e. $y=x^{-1}$.}
The above proposition motivates us to call $\G(\M)\tto \Ll(\M)$ the groupoid of partially invertible elements of $W^*$-algebra $\M$.

Now, following \cite{OS}, we  define on $\Ll(\M)$ and $\G(\M)$ structures of
complex Banach manifolds. For $p\in\Ll(\M)$  we note by $\Pi_p$ the subset
of $\Ll(\M)$ consisting of orthogonal projections $q\in\Ll(\M)$ such that the
Banach splitting \be\label{split}\M=q\M\oplus (1-p)\M\ee of $\M$ on the right
ideals $q\M$ and $(1-p)\M$ is valid. Decomposing the projection $p$ in accordance with
(\ref{split}) \be p=\sigma_p(q)-\varphi_p(q)\ee we obtain a bijective map \be
\varphi_p:\Pi_p\tilde\to(1-p)\M p\ee of $\Pi_p$ on Banach space $(1-p)\M p$ as well
as a local section \be \sigma_p:\Pi_p\to \ll^{-1}(\Pi_p)\ee of the target map, i.e.
$\ll(\sigma_p(q))=q$ for $q\in\Pi_p$. For $y\in(1-p)\M p$ we find \be
\varphi_p^{-1}(y_p)=\ll(p+y_p).\ee Using splitting (\ref{split}) for projections $p,p'\in
\Ll(\M)$ we obtain explicit expression \be\label{y_trans}
y_{p'}=\left(\varphi_{p'}\circ\varphi^{-1}_p\right)(y_p)=(b+dy_p)(a+cy_p)^{-1}\ee for the
transition map $$\varphi_{p'}\circ\varphi^{-1}_p:\varphi_p(\Pi_p\cap\Pi_{p'})\to\varphi_{p'}(\Pi_p\cap\Pi_{p'})$$ where $y_p\in(1-p)\M p$, $y_{p'}\in(1-p')\M p'$ and groupoids elements $a\in p'\M p$, $b\in (1-p')\M p$, $c\in p'\M(1-p)$, $d\in(1-p')\M(1-p)$ are defined by $\varphi_{p'}\circ\varphi^{-1}_p$ in a unique way. For the details we address to \cite{OS}. If $$ \Omega_{\tilde{p}p}:=\ll^{-1}(\Pi_{\tilde{p}})\cap \rr^{-1}(\Pi_p)\not=0$$ then one has the one-to-one map $$\psi_{\tilde{p}p}:\Omega_{\tilde{p}p}\to (1-\tilde{p})\M\tilde{p}\oplus\tilde{p}\M p\oplus(1-p)\M p$$  defined by \be\label{psipp}\psi_{\tilde{p}p}(x):=\left(\varphi_{\tilde{p}}(\ll(x)),(\sigma_{\tilde{p}}(\ll(x)))^{-1}x \sigma_{{p}}(\rr(x)),\varphi_p(\rr(x)) \right)=:(y_{\tilde p},z_{\tilde p p},y_p).\ee The image $\psi_{\tilde{p}p}(\Omega_{\tilde{p}p})$ of $\psi_{\tilde{p}p}$ is open in the Banach space $(1-\tilde{p})\M\tilde{p}\oplus\tilde{p}\M p\oplus(1-p)\M p$. For $(y_{\tilde p},z_{\tilde p p},y_p)\in\psi_{\tilde{p}p}(\Omega_{\tilde{p'}p'}\cap \Omega_{\tilde{p}p})$ one has
\be  (\psi_{\tilde{p'}p'}\circ\psi_{\tilde{p}p}^{-1})(y_{\tilde p},z_{\tilde p p},y_p):=(y_{\tilde{p'}},z_{\tilde{p'}p'},y_{p'}), \ee
 where
 \be\label{tR1}
 y_{\tilde{p'}}=(\varphi_{\tilde{p'}}\circ\varphi^{-1}_{\tilde{p}})(y_{\tilde p})=(\tilde{b}+\tilde{d}{y}_{\tilde p})(\tilde{a}+\tilde{c}{y}_{\tilde p})^{-1} \ee
 \be\label{tR2} y_{p'}=(\varphi_{p'}\circ\varphi^{-1}_{p})(y_p)=(b+dy_p)(a+cy_p)^{-1}\ee
 and
 \be\label{tR3} z_{\tilde{p'}p'}=(\tilde{p'}+y_{\tilde{p'}})^{-1}(\tilde{p}+{y}_{\tilde p})z_{\tilde p p}(p+y_p)^{-1}(p'+y_{p'}).\ee

 It follows from (\ref{y_trans}) and (\ref{tR1}-\ref{tR3}) that the atlas \be\label{atlaskrata}\left(\Pi_p,\ \varphi_p:\Pi_p\to(1-p)\M p\right),\ee where $p\in \Ll(\M)$, and the atlas  \be\label{atlasgroup}\left(\Omega_{\tilde{p}p},\psi_{\tilde{p}p}:\Omega_{\tilde{p}p}\to (1-\tilde{p})\M\tilde{p}\oplus\tilde{p}\M p\oplus(1-p)\M p\right),\ee where $\tilde p,p\in\Ll(\M)$, define structures of complex Banach manifolds on the lattice $\Ll(\M)$ and on the groupoid $\G(\M)$, respectively. In \cite{OS} it was proved that $\G(\M)\tto\Ll(\M)$ is the Banach-Lie groupoid with respect to the complex Banach manifold structures defined by the atlases (\ref{atlaskrata}) and (\ref{atlasgroup}).

 It is reasonable to note here that if $p'\nsim p$ then $\Pi_{p'}\cap\Pi_p=\emptyset$. However the equivalence $p'\sim p$ of the orthogonal projections does not imply $\Pi_{p'}\cap\Pi_p\not=\emptyset$. For example if $\M$ is an infinite $W^*$-algebra then there exists such a projection that $\mathbf{1}\not=p\sim \mathbf{1}$. For this projection we have $\Pi_\mathbf{1}=\{\mathbf{1}\}$ and thus $\Pi_\mathbf{1}\cap\Pi_p=\emptyset$.

 We also note that $\Omega_{p'p}\not=\emptyset$ if and only if $p'\sim p$. So, von Neumann equivalence of projections $p'$ and $p$ means that they belongs to the same orbit  of the canonical action $\G(\M)\times\Ll(\M)\to\Ll(\M)$ of the groupoid $\G(\M)$ on the lattice $\Ll(\M)$.

  Let $\Ll_{p_0}(\M)\subset\Ll(\M)$ be the orbit of this action generated by $p_0\in\Ll(\M)$ i.e.
 \be\label{orbit} \Ll_{p_0}(\M):=\{\ll(x):\quad x\in\G(\M),\quad \rr(x)=p_0\}.\ee

 By $\G_{p_0}(\M)\tto\Ll_{p_0}(\M)$ we denote the transitive subgroupoid of $\G(\M)\tto\Ll(\M)$, where \be\label{Gp0}\G_{p_0}(\M):=\ll^{-1}(\Ll_{p_0}(\M))\cap \rr^{-1}(\Ll_{p_0}(\M)).\ee
 Since for $\tilde p\nsim p$ one has $\Pi_{p'}\cap\Pi_p=\emptyset$ and $\Omega_{p'p}=\emptyset$ then the locally trivial  groupoid $\G_{p_0}(\M)\tto\Ll_{p_0}(\M)$ is a Banach-Lie subgroupoid of the groupoid $\G(\M)\tto\Ll(\M)$.  One shows the local triviality of $\G_{p_0}(\M)\tto\Ll_{p_0}(\M)$ expressing $r$ and $l$ in coordinates $(y_{\tilde p}, z_{\tilde p p},y_p)$.

 Note that $\G(\M)\tto\Ll(\M)$ is a union of the transitive  subgroupoids $\G_{p_0}(\M)\tto\Ll_{p_0}(\M)$ and this decomposition is consistent with its complex Banach manifold structure. Summing up we can reduce the investigations of $\G(\M)\tto\Ll(\M)$ to the investigations of the transitive subgroupoids (\ref{Gp0}).

 According to \cite{OS} let us define the involution map $J:\G(\M)\to\G(\M)$ on the groupoid $\G(\M)\tto\Ll(\M)$ in the following way
 \be\label{J} J(x):=(x^*)^{-1}.\ee
 From the polar decomposition (\ref{polar}) of $x\in \G(\M)$ we find that $J(x)=x$ iff $|x|=r(x)$, i.e. fixed points of $J$ form the subset $\U(\M)\subset \G(\M)$ of partial isometries which is invariant with respect to groupoid operations. In particular for $u\in \U(\M)$ we have
 \be u^{-1}=u^*,\quad r(u)=u^*u,\quad l(u)=uu^*\ee
 and the base set as well as the objects inclusion map $\varepsilon:\Ll(\M)\to \U(\M)$ are the same as for $\G(\M)\tto\Ll(\M)$. So, $\U(\M)\tto\Ll(\M)$ is a wide subgroupoid of the groupoid of partially invertible elements.

 We note that the involution map $J$ preserves $\Omega_{\tilde pp}$ and in the coordinates (\ref{psipp}) it is given by
 \be \left(\psi_{\tilde{p}p}\circ J\circ \psi_{\tilde{p}p}^{-1}\right)\left( y_{\tilde p},z_{\tilde pp},y_p\right)=\ee
 $$=\left( y_{\tilde p},((\tilde p+y_{\tilde p})^*(\tilde p+y_{\tilde p}))^{-1}z^*_{\tilde pp}(p+y_p)^*(p+y_p),y_p\right).$$

  In what follows  subsequently we will denote the transitive groupoid $\left(\U(\M)\cap\G_{p_0}(\M)\right)\tto\Ll(\M)$ by $\U_{p_0}(\M)\tto\Ll_{p_0}(\M)$. The smooth Banach manifold structure of $\U_{p_0}(\M)\tto\Ll_{p_0}(\M)$ we will describe in Section 3.

 \section{Algebroids of the groupoids $\G_{p_0}(\M)\tto\Ll_{p_0}(\M)$ and $\U_{p_0}(\M)\tto\Ll_{p_0}(\M)$}

 At first we will discuss several questions concerning the structure of $\G_{p_0}(\M)\tto\Ll_{p_0}(\M)$. The  geometric constructions investigated here turn out to be useful for the description of the Banach-Lie algebroid $\A_{p_0}(\M)$ of $\G_{p_0}(\M)\tto\Ll_{p_0}(\M)$. Let us introduce the following notations. By $G_0$ we denote the group of invertible elements of $W^*$-subalgebra $p_0\M p_0\subset\M$. The intersection $\G_{p_0}(\M)\cap \M p_0$ of $\G_{p_0}(\M)$ with the left $W^*$-ideal $\M p_0$ we denote by $P_0$. We indicate the following statement.
 \begin{prop}\label{21}
 \ben[(i)]
 \item The group $G_0$ is an open subset of the Banach space $p_0\M p_0$. So, $G_0$ is a Banach-Lie group whose Lie algebra is $p_0\M p_0$.
     \item The subset $P_0\subset \M p_0$  is open in the Banach space $\M p_0$. Thus the tangent vector bundle $TP_0$  can be identified with the trivial vector bundle $\M p_0\times P_0$.
 \item One has a free right action of $G_0$ on $P_0\times P_0$ given by
 \be\label{G0action}  P_0\times P_0\times G_0\ni (\eta,\xi,g)\ \mapsto\ (\eta g,\xi g)\in P_0\times P_0.\ee\een
 \end{prop}
 The left support map $\ll:P_0\to \Ll_{p_0}(\M)$ is a surjective submersion of Banach manifolds which is invariant with respect to the right action
  \be\label{actionGP} P_0\times G_0\ni (\eta,g)\ \mapsto\ \eta g\in P_0\ee
of $G_0$ on $P_0$,  i.e. for $g\in G_0$ and $\eta\in P_0$ one has $\ll(\eta g)=\ll(\eta)$. We conclude from the above that \be\label{baza}\Ll_{p_0}(\M)\cong {P_0}/{G_0}\ee  and  $P_0\left(\Ll_{p_0}(\M), \ G_0,\ l\right)$ defines a principal bundle structure on $P_0$.

 According to item $(ii)$ of the  Proposition \ref{21} one can identify the element $\eta \in P_0\hookrightarrow\M {p_0}$ with its coordinate in Banach space $\M {p_0}$. However the inclusion map $\iota: P_0\hookrightarrow \M p_0$ is not consistent with the principal bundle  structure of $P_0$. Therefore  we will use further the atlas given by one-to-one maps
 $ \psi_p:\ll^{-1}(\Pi_p)\to (1-p)\M p\oplus p\M p_0$, for $p\in\Ll_{p_0}(\M)$,
defined as follows
 \be\label{atlasP}\psi_p(\eta):=\left(\varphi_p(\ll(\eta)), \sigma_p(\ll(\eta))^{-1}\eta\right)=\ee
 $$=\left(\varphi_p(\eta\eta^{-1}), (p+\varphi_p(\eta\eta^{-1}))^{-1}\eta\right)=:(y_p,z_p).$$
 The transition maps $$\psi _{p'}\circ \psi_p^{-1}:\psi_p(\ll^{-1}(\Pi_p\cap\Pi_{p'}))\to \psi_{p'}(\ll^{-1}(\Pi_p\cap\Pi_{p'}))$$ for this atlas are given by
 \be y_{p'}=(b+dy_p)(a+cy_p)^{-1},\ee
 \be z_{p'}=(p'+y_{p'})^{-1}(p+y_p)z_p.\ee

 The atlas (\ref{atlasP}) is consistent with the groupoid atlas given in (\ref{atlasgroup}). Namely, we obtain the chart (\ref{atlasP}) as $\psi_p=\psi_{pp_0}$ taking  $\tilde p=p$, $\ p=p_0$ and $y_{p_0}=0$ in (\ref{psipp}).  Note here that the set $\psi_p(\ll^{-1}(\Pi_p))$ is open in $(1-p)\M p\oplus p\M p_0$. Note also that the map $\psi_p^{-1}$ inverse  to $\psi_p$ has the form
 \be\label{eta} \eta=\psi_p^{-1}(y_p,z_p)=(p+y_p)z_p.\ee

 Inverting (\ref{eta}) we find that
 \be\label{yetazeta} y_p=(1-p)\eta(p\eta)^{-1},\qquad z_p=p\eta.\ee

 In order to see that the dependences (\ref{yetazeta}) are complex analytic we assume that  $\norm{\lambda_{pp_0}-p\eta}\leqslant\norm{\lambda^{-1}_{pp_0}}^{-1}$. Then applying the Neumann formula to $\lambda^{-1}_{pp_0}p\eta\in G_0$ we obtain
 \be (p\eta)^{-1}=\sum_{n=0}^\infty \left[\lambda^{-1}_{pp_0}(\lambda_{pp_0}-p\eta)\right]^n\lambda_{pp_0}.\ee
 So, the chart $(P_0, \iota:P_0\hookrightarrow \M p_0)$ belongs to the maximal atlas generated by the charts $(l^{-1}(\Pi_p),\psi_p)$, $\ p\in \Ll_{p_0}(\M).$

 Now for any projection $p\in \Ll_{p_0}(\M)$ let us fix a groupoid element $\lambda_{pp_0}\in p\M p_0\cap\G_{p_0}(\M)$. The maps
 $\Lambda_p:\ll^{-1}(\Pi_p)\to \Pi_p\times G_0$ defined by
 \be\label{maps} \Lambda_p(\eta):=\left(\ll(\eta), \lambda^{-1}_{pp_0}\sigma_p(\ll(\eta))^{-1}\eta\right)\ee
 give a local trivialization of the principal bundle $P_0\left(\Ll_{p_0}(\M), \ G_0,\ l\right)$. The transition cocycle $ \mathbf{g}_{p'p}:\Pi_{p'}\cap\Pi_p\to G_0$ corresponding to (\ref{maps}) has the form
 \be \mathbf{g}_{p'p}(q)=\lambda^{-1}_{p'p_0}\sigma_{p'}^{-1}(q)\sigma_p(q)\lambda_{pp_0}.\ee
 Another choice of $\{\lambda_{pp_0}\}$ gives a cocycle equivalent to the cocycle $\{\mathbf{g}_{p'p}\}$. One has the following relations
 \be\left(\Lambda_p\circ \psi^{-1}_p\right)(y_p,z_p)=(\varphi ^{-1}_p(y_p),\lambda^{-1}_{pp_0}z_p).\ee
 So one can use the coordinates $(y_p,\zeta_p)\in (1-p)\M p\oplus p_0\M p_0$, where
 \be\zeta_p:=\lambda^{-1}_{pp_0}z_p\ee
  instead of the coordinates $(y_p,z_p)\in (1-p)\M p\oplus p\M p_0$. The above shows that the Banach manifold structure of $P_0\left(\Ll_{p_0}(\M), \ G_0,\ l\right)$ can be modeled on the Banach spaces $(1-p)\M p\oplus p_0\M p_0$, where $p\in \Ll_{p_0}(\M)$.

  In the general  theory of Lie groupoids an important role is played by the gauge groupoid $\frac{P\times P}{G}\rightrightarrows  M$ associated in a canonical way to a principal bundle $P(M, G, \pi)$, e.g. see \cite{mac}. In our case  we obtain the groupoid $\frac{P_0\times P_0}{G_0}\tto {P_0}/{G_0}$, where for $\langle\eta,\xi\rangle:=\{(\eta g,\xi g);\ \ g\in G_0\}\in \frac{P_0\times P_0}{G_0}$ and $[\eta]:=\{\eta g:\ g\in G_0\}\in {P_0}/{G_0}$ one defines
  \be \begin{array}{l}
  \Ss(\langle\eta,\xi\rangle):=[\xi]\\
  \Tt(\langle\eta,\xi\rangle):=[\eta]\\
  \langle\eta,\xi\rangle^{-1}:=\langle\xi,\eta\rangle\\
  \varepsilon([\eta]):=\langle\eta,\eta\rangle\\
  \langle\eta,\xi\rangle\cdot\langle\lambda,\kappa\rangle=\langle\eta,\kappa g\rangle,\end{array}\ee
 where $g\in G_0$ is given by $\xi=\lambda g$.
  \begin{prop}\label{prop22} One has the following isomorphism
   \unitlength=5mm \be\label{gaugethm}\begin{picture}(11,4.6)
    \put(1,4){\makebox(0,0){$\frac{P_0\times P_0}{G_0}$}}
    \put(8,4){\makebox(0,0){$\G_{p_0}(\M)$}}
    \put(1,-1){\makebox(0,0){$P_0/G_0$}}
    \put(8,-1){\makebox(0,0){$\Ll_{p_0}(\M)$}}
    \put(1.2,3){\vector(0,-1){3}}
    \put(0.7,3){\vector(0,-1){3}}
    \put(8.2,3){\vector(0,-1){3}}
    \put(7.7,3){\vector(0,-1){3}}
    \put(3,4){\vector(1,0){3}}
    \put(2.7,-1){\vector(1,0){3.7}}
    \put(0.1,1.4){\makebox(0,0){$\Ss$}}
    \put(2.2,1.4){\makebox(0,0){$\Tt$}}
    \put(9.1,1.4){\makebox(0,0){$l$}}
    \put(6.8,1.4){\makebox(0,0){$r$}}
    \put(4.5,4.5){\makebox(0,0){$\phi$}}
    \put(4.5,-0.5){\makebox(0,0){$\varphi $}}
    \end{picture},\ee

    \bigskip

  of Banach-Lie groupoids, where
  \be\label{I}
 \phi:\frac{P_0\times P_0}{G_0}\ \ni\langle{\eta},{\xi}\rangle\mapsto\eta\xi^{-1}\  \in\ \G_{p_0}(\M),\ee
 and
 \be \label{II}\varphi:P_0/G_0\ \ni\ [\eta]\ \mapsto\ \eta\eta^{-1}\in \Ll_{p_0}(\M).\ee
 \end{prop}

 Note here that the coordinates $(y_{\tilde p}, z_{\tilde p p}=z_{\tilde p}z_p^{-1}, y_p) $ of $\eta\xi^{-1}\in\G_{p_0}(\M)$ defined in (\ref{psipp}) are invariant with respect to the right action (\ref{G0action}) of the group $G_0$ on the $P_0\times P_0$.

 Now let us mention several facts concerning the group $(TG_0,\bullet)$ which is the tangent group to the group $G_0$. The product of $X_g\in T_g G_0$ and $Y_h\in T_h G_0$ is given by
 \be\label{tangent_prod} X_g\bullet Y_h=TR_h(g)X_g+TL_g(h)Y_h.\ee
 The  space $T_eG_0\cong p_0\M p_0$ tangent to $G_0$ at the unit element $e=p_0$ is a normal subgroup of $TG_0$.
 One has the  decomposition $TG_0=T_eG_0\bullet G_0$ of $TG_0$, where we identify $G_0$ with the zero section of $TG_0$. Hence we have the following group isomorphisms
 \be\label{TG0} TG_0= G_0\ltimes_{Ad} T_e G_0\quad {\rm and} \quad G_0\cong TG_0/T_e G_0.\ee
 Using (\ref{TG0}) we obtain from (\ref{tangent_prod}) that the product of $(g,x), (h,y) \in G_0\ltimes_{Ad} p_0\M p_0$ is expressed by the formula
 \be (g,x)\bullet (h,y)=(gh,x+gyg^{-1}).\ee
 The right action of the tangent group $TG_0$ on the tangent vector bundle $TP_0\cong \M p_0\times P_0$ defined by the action (\ref{actionGP})  is given by
 \be\label{tangent_action} (\vartheta,\eta)\ast(g,x)=(\vartheta g+\eta x g, \eta g),\ee
 where $(\vartheta,\eta)\in \M p_0\times P_0$.

 The vector subspace $T_\eta^V P_0\subset T_\eta P_0$ tangent to the orbit $[\eta]$ at $\eta\in[\eta]\subset P_0$ is equal to
 \be (\eta p_0\M p_0,\eta)=\{(0,\eta)\ast(p_0,x):\quad (p_0,x)\in T_{p_0}G_0\}.\ee
 Thus we have the following isomorphisms of vector bundles
 \be\label{Tv} T^VP_0\cong (\{0\}\times P_0)\ast T_{p_0}G_0\ee
 \be\label{Tviso} TP_0/T^VP_0\cong TP_0/T_{p_0}G_0.\ee
 The action of $G_0$ on the vector bundles over $P_0$ included into the exact sequence
 \be\label{TPexact} 0\to T^VP_0\hookrightarrow TP_0\rightarrowtail TP_0/T^VP_0\to 0\ee commutes with the morphisms of (\ref{TPexact}). Thus we obtain the short exact sequence
 \be\label{TPexactiloraz} 0\to T^VP_0/G_0\stackrel{\iota}\to TP_0/G_0\stackrel{\mathbf{a}}\to TP_0/TG_0\to 0\ee of the factor vector bundles over ${P_0}/{G_0}$. In order to obtain (\ref{TPexactiloraz}) we used the isomorphism (\ref{Tviso}) and the decomposition $TG_0=T_eG_0\bullet G_0$.

In  such a way we get the Atiyah sequence
\be\label{atyah} 0\to p_0\M p_0\times_{Ad_{G_0}} P_0\stackrel{\iota}{\hookrightarrow} TP_0/ {G_0}\stackrel{\mathbf{a}}{\rightarrow} T(P_0/G_0)\to 0\ee
of the principal bundle $P_0\left(\Ll_{p_0}(\M), \ G_0,\ l\right)$,  where $p_0\M p_0$ is the Lie algebra of the Banach-Lie group $G_0$, $TP_0/G_0$ is the Banach Lie algebroid of the groupoid $\frac{P_0\times P_0}{G_0}\tto {P_0}/{G_0}$ and  \be\label{bazizo}T({P_0}/{G_0})\cong {TP_0}/{TG_0}\ee is the tangent bundle of ${P_0}/{G_0}$. The vector bundle monomorphism $\iota$ and bundle epimorphism $\mathbf{a}$ are defined by the quotient of (\ref{TPexact}). It follows from the Proposition \ref{prop22} that
\be TP_0/G_0 \cong  \A _{p_0}(\M)\ee
and
\be T(P_0/G_0)\cong T\Ll_{p_0}(\M),\ee
where $\A _{p_0}(\M)\stackrel{\pi}{\rightarrow} \Ll_{p_0}(\M)$ is the Banach-Lie algebroid of $\G_{p_0}(\M)\tto\Ll_{p_0}(\M)$ and $\mathbf{a}$ by virtue of (\ref{atyah}) is the anchor map for $\A _{p_0}(\M)$.

 In order to identify the Lie bracket structure of  space $\Gamma^\infty(TP_0/G_0)$ of smooth real  sections of the vector bundle $TP_0/G_0\stackrel{l}{\longrightarrow}\Ll_{p_0}(\M)$, we notice that these sections can be considered as $G_0$-invariant real vector fields
 $\mathfrak{X}\in\Gamma^\infty_{G_0}(TP_0)$ on $P_0$, i.e.
 \be\label{X} \mathfrak{X}(\eta,\eta^*)=\langle\frac{\partial}{\partial\eta},\vartheta(\eta,\eta^*)\rangle
 +\langle\frac{\partial}{\partial\eta^*},\vartheta(\eta,\eta^*)\rangle\ee
 and
 \be\label{varth} \vartheta(\eta g,(\eta g)^*)=\vartheta(\eta,\eta^*)g,\ee where $\vartheta:P_0\to\M p_0$ is a smooth map. We use in (\ref{X}) the complex coordinates $(\eta,\eta^*)\in  \M p_0\oplus p_0\M$. Let us note that in the real coordinates $\chi:=\frac{1}{2}(\eta+\eta^*), \ \nu:=\frac{1}{2i}(\eta-\eta^*)\in (\M p_0\oplus p_0\M)^h$, where $(\M p_0\oplus p_0\M)^h$ is the hermitian part of complex Banach space $\M p_0\oplus p_0\M$, the vector field (\ref{X}) is given by
 \be\label{Xreal} \mathfrak{X}(\chi,\nu)=\langle \frac{\partial}{\partial\chi},\alpha(\chi,\nu)\rangle+\langle \frac{\partial}{\partial\nu},\beta(\chi,\nu)\rangle,\ee
 where $\alpha=\frac{1}{2i}(\vartheta+\vartheta^*)$ and $\beta=\frac{1}{2i}(\vartheta-\vartheta^*)$.

 We explain the notation in (\ref{X}) and (\ref{Xreal}). Let $\gamma_t:P_0\to P_0, \ t\in \mathbb{R}$ be a (local) one-parameter group of automophisms of the principal bundle $P_0(\Ll_{p_0}(\M),G_0,l)$ which is tangent to the real vector field $\mathfrak{X}$, i.e. in the complex coordinates $(\eta, \eta^*)\in \M p_0\oplus p_0\M$ it satisfies
 \be \gamma_t(\eta g)=\gamma_t(\eta)g\ee
 and
 \be \frac{d}{dt}\gamma_t(\eta)|_{t=0}=\frac{d}{dt}\gamma(t,\eta,\eta^*)|_{t=0}=\vartheta(\eta,\eta^*).\ee
 Then for real valued smooth function $f\in \mathcal{C}^\infty(P_0,\mathbb{R})$ one has
 \be\label{Xf} ( \mathfrak{X}f)(\eta, \eta^*)=\frac{d}{dt}f(\gamma_t(\eta), \gamma_t(\eta)^*)|_{t=0}=\ee
 $$=\langle \frac{\partial f}{\partial \eta}(\eta, \eta^*),  \frac{d}{dt}\gamma_t(\eta)|_{t=0}\rangle+
 \langle \frac{\partial f}{\partial \eta^*}(\eta, \eta^*),  \frac{d}{dt}\gamma_t(\eta)|_{t=0}^*\rangle=$$
 $$=\langle \frac{\partial f}{\partial \eta}(\eta, \eta^*),  \vartheta(\eta,\eta^*)\rangle+
 \langle \frac{\partial f}{\partial \eta^*}(\eta, \eta^*),  \vartheta(\eta,\eta^*)^*\rangle=$$
 $$=\langle \frac{\partial f}{\partial\chi}(\chi,\nu),\alpha(\chi,\nu)\rangle+\langle \frac{\partial f}{\partial\nu}(\chi,\nu),\beta(\chi,\nu)\rangle.$$
 Note here that
 $\frac{\partial f}{\partial \eta}(\eta, \eta^*)\in (\M p_0)^*$,
 $\  \vartheta(\eta,\eta^*)\in \M p_0$,
 $\ \frac{\partial f}{\partial \eta^*}(\eta, \eta^*)\in (p_0\M)^*$,
 $\  \vartheta(\eta,\eta^*)^*\in p_0\M$,
 $\alpha(\chi,\nu),$ $ \beta(\chi,\nu)\in \M p_0\oplus p_0\M$ and $\frac{\partial f}{\partial\chi}(\chi,\nu),$ $\frac{\partial f}{\partial\nu}(\chi,\nu)\in ((\M p_0\oplus p_0\M)^h)^*$. So, the paring $\langle\cdot , \cdot\rangle$ on the right hand side of (\ref{Xf}) is correctly defined. Removing the function $f$ from (\ref{Xf}) we obtain (\ref{X}).
 
 We note here that $(\M p_0)^*$, $\ (p_0\M)^*$ and $((\M p_0\oplus p_0\M)^h)^*$ are Banach spaces dual to the Banach subspaces $\M p_0$, $p_0\M$ and $(\M p_0\oplus p_0\M)^h$, respectively. The analogous convention we will also use  in the subsequent.

 In this notation the Lie bracket of $\mathfrak{X}_1,\ \mathfrak{X}_2\in\Gamma^\infty_{G_0}(TP_0)$ is given by
 $$ [\mathfrak{X}_1,\mathfrak{X}_2]=\left\langle\frac{\partial }{\partial \eta},\langle\frac{\partial \vartheta_2}{\partial \eta}, \vartheta_1\rangle-\langle\frac{\partial \vartheta_1}{\partial \eta}, \vartheta_2\rangle+
 \langle\frac{\partial \vartheta_2}{\partial \eta^*}, \vartheta_1^*\rangle-\langle\frac{\partial \vartheta_1}{\partial \eta^*}, \vartheta_2^*\rangle\right\rangle+$$
\be\label{bracket} +\left\langle\frac{\partial }{\partial \eta^*},\langle\frac{\partial \vartheta_2^*}{\partial \eta}, \vartheta_1\rangle-\langle\frac{\partial \vartheta_1^*}{\partial \eta}, \vartheta_2\rangle+
 \langle\frac{\partial \vartheta_2^*}{\partial \eta^*}, \vartheta_1^*\rangle-\langle\frac{\partial \vartheta_1^*}{\partial \eta^*}, \vartheta_2^*\rangle\right\rangle\ee

It is easy to see that $[\mathfrak{X}_1,\mathfrak{X}_2]\in \Gamma^\infty_{G_0}(TP_0)\cong \Gamma^\infty(TP_0/G_0)$.

 The vector field (\ref{X}) is associated to the one-parameter group $L_t\circ L_s=L_{t+s}$ of the left translation $L_t:\G_{p_0}(\M)\to\G_{p_0}(\M)$ of the groupoid $\G_{p_0}(\M)\tto\Ll_{p_0}(\M)$ in the following way. From the defining  properties
\be\begin{array}{rl}
 (i)& L_t(xy)=L_t(x)y\\
(ii)& r\circ L_t=r\\
(iii)& l\circ L_t=\lambda_t\circ l\end{array} \ee
 of $L_t$, where $\rr(x)=\ll(y)$ and $\lambda_t:\Ll_{p_0}(\M)\to\Ll_{p_0}(\M)$ is a one-parameter group of smooth transformations of $\Ll_{p_0}(\M)$, we see that $\{L_t\}_{t\in\mathbb{R}}$ as well as $\{\lambda_t\}_{t\in\mathbb{R}}$ are defined by a one-parameter group of automorphisms  $\gamma_t:P_0\to P_0$ of the principal bundle $P_0\left(\Ll_{p_0}(\M), \ G_0,\ l\right)$, i.e.
 \be L_t(\eta\xi^{-1})=\gamma_t(\eta)\xi^{-1}.\ee

 The algebroid properties of the bracket (\ref{bracket}) are better seen in the coordinates $(y_p,z_p,y_p^*,z_p^*)\in (1-p)\M p\oplus p\M p_0\oplus p\M(1-p)\oplus p_0\M p$. For this reason let us notice that in these coordinates   the flow $\gamma_t:P_0\to P_0$  has the form
 \be\label{yp} \mathbf{y}_p(t,{y}_p, {y}_p^*):= \left(\varphi_p\circ \lambda_t\circ\varphi_p^{-1}\right)(y_p)=y_p(t)\ee
 and
 \be\label{z} \mathbf{z}_p(t,y_p,y_p^*,z_p):=
 \left(p+\mathbf{y}_p(t,{y}_p, {y}_p^*)\right)^{-1}\gamma_t\left(p+y_p)\lambda_{pp_0}\right)\lambda_{pp_0}^{-1}z_p=\qquad\qquad\ee
 $$\qquad= \left(p+\mathbf{y}_p(t,{y}_p, {y}_p^*)\right)^{-1}\gamma(t,(p+y_p)\lambda_{pp_0},\lambda_{pp_0}^*(p+y_p)^*)\lambda_{pp_0}^{-1}z_p=z_p(t),$$
 where $t\in]-\varepsilon,\varepsilon[$.  \\ The equality (\ref{z}) we obtain from
 \be (p+y_p(t))z_p(t)=\gamma_t(\eta)=\gamma_t\left((p+y_p)\lambda_{pp_0}\lambda_{pp_0}^{-1}z_p\right)=\ee
 $$=\gamma_t\left((p+y_p)\lambda_{pp_0}\right)\lambda_{pp_0}^{-1}z_p.$$
 From (\ref{yp}) and (\ref{z}) it follows that
  $$\mathfrak{X}(y_p,z_p,y_p^*,z_p^*)=\langle \frac{\partial }{\partial y_p},a_p(y_p, y_p^*)\rangle+\langle \frac{\partial }{\partial z_p},b_p(y_p, y_p^*)z_p\rangle+$$
 \be\label{frakX}+\langle \frac{\partial }{\partial y_p^*},a_p(y_p, y_p^*)^*\rangle+\langle \frac{\partial }{\partial z_p^*},z_p^* b_p(y_p, y_p^*)^*\rangle,\ee
 where $a_p:(1-p)\M p\to(1-p)\M p$ and $b_p:(1-p)\M p\to p\M p$ are smooth maps defined by
 \be \begin{array}{l} a_p(y_p,y_p^*):=\frac{d}{dt}\mathbf{y}_p(t,y_p,y_p^*)|_{t=0}\\
 b_p(y_p, y_p^*)z_p:=\frac{d}{dt}\mathbf{z}_p(t,y_p,y_p^*,z_p)|_{t=0}\end{array}\ee
 The symbols $\frac{\partial }{\partial y_p}$, $\frac{\partial }{\partial y_p^*}$, $\frac{\partial }{\partial z_p}$  and $\frac{\partial }{\partial z_p^*}$ in (\ref{frakX}) denote the corresponding partial derivatives.
 Rewriting (\ref{bracket}) in the coordinates $(y_p,z_p,y_p^*,z_p^*)$ we find that the bracket $[\mathfrak{X}_1,\mathfrak{X}_2]$ is given by
\be\label{nawias}[\mathfrak{X_1},\mathfrak{X_2}]
=\left\langle \frac{\partial}{\partial y_p},a_p\right\rangle+\left\langle \frac{\partial}{\partial y_p^*},a_p^*\right\rangle+\left\langle \frac{\partial}{\partial z_p},b_pz_p\right\rangle+\left\langle \frac{\partial}{\partial z_p^*},z_p^*b_p^*\right\rangle,\ee
where
\be\label{ap} a_p:=\langle\frac{\partial a_{2p}}{\partial y_p},a_{1p}\rangle-\langle\frac{\partial a_{1p}}{\partial y_p},a_{2p}\rangle+\langle\frac{\partial a_{2p}}{\partial y_p^*},a_{1p}^*\rangle-\langle\frac{\partial a_{1p}}{\partial y_p^*},a_{2p}^*\rangle\ee
and
\be\label{bp} b_p:=\langle\frac{\partial b_{2p}}{\partial y_p},a_{1p}\rangle-\langle\frac{\partial b_{1p}}{\partial y_p},a_{2p}\rangle+\langle\frac{\partial b_{2p}}{\partial y_p^*},a_{1p}^*\rangle-\langle\frac{\partial b_{1p}}{\partial y_p^*},a_{2p}^*\rangle+[b_{2p},b_{1p}].\ee

 The left support map $\ll:\ll^{-1}(\Pi_p)\to\Pi_p$ in the coordinates $(y_p,z_p, y_p^*,z_p^*)$ is  the projection of $(y_p,z_p)$ on the first component, i.e.  $(\varphi_p\circ\ll\circ \psi_p^{-1})(y_p,z_p)= y_p$ . The anchor map $\mathbf{a}$ is the tangent map
$$\mathbf{a}:=T\ll:\Gamma^\infty(TP_0/G_0)\to \Gamma^\infty(T\Ll_{p_0}(\M))$$ of $\ll$ and  maps
$\mathfrak{X}$ on
\be\label{anchorT}\mathbf{a}(\mathfrak{X})(y_p,y_p^*)=\langle \frac{\partial }{\partial y_p},a_p(y_p, y_p^*)\rangle+\langle \frac{\partial }{\partial y_p^*},a_p(y_p,y_p^*)^*\rangle.\ee
So, the property
\be\mathbf{a}([\mathfrak{X}_1,\mathfrak{X}_2])=[\mathbf{a}(\mathfrak{X}_1),\mathbf{a}(\mathfrak{X}_2)]\ee
follows from (\ref{nawias}), (\ref{anchorT}). The Jacobi identity and the Leibniz property
\be [\mathfrak{X}_1,f \mathfrak{X}_2]=f [\mathfrak{X}_1,\mathfrak{X}_2]+\mathbf{a}(\mathfrak{X}_1)(f)\mathfrak{X}_2,\ee
where $f\in\mathcal{C}^\infty({P_0}/{G_0})=\mathcal{C}_{G_0}^\infty(P_0)$, can be easily seen when written in the coordinates $(y_p,z_p,y_p^*,z_p^*)$.

The vertical $G_0$-invariant vector field $\mathfrak{X}\in \Gamma^\infty_{G_0}(T^VP_0)$ is given by
\be\label{XV} \mathfrak{X}^V(y_p,z_p,y_p^*,z_p^*)=\langle\frac{\partial}{\partial z_p},b_p(y_p,y_p^*)z_p\rangle +
\langle\frac{\partial}{\partial z_p^*},z_p^*b_p(y_p,y_p^*)^*\rangle.\ee
It follows from (\ref{ap}) that the space of $G_0$-invariant vector fields $\Gamma^\infty_{G_0}(T^VP_0)$ as an ideal in the algebroid $\Gamma^\infty_{G_0}(TP_0)$, what is in agreement with general theory of algebroids.

Let us also mention that the vector subspace $\Gamma^{hol}_{G_0} (TP_0)\subset\Gamma^\infty_{G_0} (TP_0)$ of holomorphic vector fields, i.e. $\mathfrak{X}\in \Gamma^{hol}_{G_0} (TP_0)$ if and only if $\frac{\partial\vartheta}{\partial\eta^*}\equiv 0$, is a subalgebroid of the algebroid $\Gamma^\infty_{G_0} (TP_0)$.

Since
\be \frac{\partial }{\partial y_p}=\frac{\partial \eta}{\partial y_p}\frac{\partial }{\partial \eta}=z_p\frac{\partial }{\partial \eta}\ee
and
\be  \frac{\partial }{\partial z_p}=\frac{\partial \eta}{\partial z_p}\frac{\partial }{\partial \eta}=(p+y_p)\frac{\partial }{\partial \eta}\ee
we find that
\be\label{wzor} \vartheta=a_p(y_p,y_p^*)z_p+(p+y_p)b_p(y_p,y_p^*)z_p\ee
and
\be\label{ab} a_p=(1-\eta(p\eta)^{-1})\vartheta(p\eta)^{-1},\ee
\be\label{ba} b_p=p\vartheta (p\eta)^{-1}.\ee

Let us remark here that for the proof of the equalities (\ref{wzor}), (\ref{ab}) and (\ref{ba}) it is important to keep the proper order of the corresponding  $\M$-valued (operator valued) functions and remember that $(p\eta)^{-1}\not=\eta^{-1}p$.

Given the vector field $\mathfrak{X}\in\Gamma^\infty_{G_0} (TP_0)$, see (\ref{X}), the corresponding section ${X}\in\Gamma^\infty \A_{p_0}(\M)$ of the algebroid $\A_{p_0}(\M)$ is defined by
\be\label{Xeta} X(\ll(\eta))=\vartheta(\eta,\eta^*)\eta^{-1}.\ee

From (\ref{bracket}) and (\ref{Xeta}) we find that the algebroid bracket $[\cdot,\cdot]_{\A}$ of $X_1, X_2\in \Gamma^\infty\A_{p_0}(\M)$ has the following form
\be\label{nawias3}[X_1,X_2]_{\A}=[X_2,X_1]+\ee
$$+\langle\frac{\partial X_2}{\partial\eta}, X_1\eta\rangle-\langle\frac{\partial X_1}{\partial\eta}, X_2\eta\rangle+\langle\frac{\partial X_2}{\partial\eta^*},\eta^* X_1^*\rangle-
\langle\frac{\partial X_1}{\partial\eta^*},\eta^* X_2^*\rangle.$$

Using the chart $(\Pi_p,\varphi_p:\Pi_p\to (1-p)\M p)$ one expresses the section $X\in \Gamma^\infty\A_{p_0}(\M)$ locally by the map
\be X_p:\varphi_p(\Pi_p)\ \ni\ y_p \mapsto X_p(y_p,y_p^*)\ \in \ \pi^{-1}(\Pi_p)\subset\M,\ee
where $\pi:\A_{p_0}(\M)\to \Ll_{p_0}(\M)$ is the bundle map of the algebroid. So, rewriting (\ref{nawias3}) in the coordinates $(y_p,y_p^*)$ we obtain
\be\label{nawias5} [X_{1p},X_{2p}]_\A=[X_{2p},X_{1p}]+\ee
$$+\langle\frac{\partial X_2}{\partial y_p},(X_1-(p+y_p)X_1)(p+y_p)\rangle -
\langle\frac{\partial X_1}{\partial y_p},(X_2-(p+y_p)X_2)(p+y_p)\rangle+$$
$$+\langle\frac{\partial X_2}{\partial y_p^*},[(X_1-(p+y_p)X_1)(p+y_p)]^*\rangle-
\langle\frac{\partial X_1}{\partial y_p^*},[(X_2-(p+y_p)X_2)(p+y_p)]^*\rangle.$$

\bigskip

At the end of this section we describe the Banach-Lie structure of the subgroupoid of partial isometries $\U_{p_0}(\M)\tto\Ll_{p_0}(\M)$. For this reason let us consider $U_0$-principal bundle $P_0^u(\Ll_{p_0}(\M),U_0,l)$, where $U_0:=U(p_0\M p_0)$ is the group of unitary elements of the $W^*$-subalgebra $p_0\M p_0$ and the total bundle space $P_0^u\subset P_0$ consisting of partial isometries of $P_0$, i.e.
\be\label{PU} \eta\in P_0^u\quad {\rm iff} \quad \eta^*\eta=p_0.\ee
Hence, by virtue of the real Banach spaces splitting
\be\label{splitting} p_0\M p_0=i p_0\M ^hp_0 \oplus  p_0\M^h p_0\ee
it follows that (\ref{PU}) that $P_0^u(\Ll_{p_0}(\M),U_0,l)$ is the principal subbundle of the principal bundle $P_0(\Ll_{p_0}(\M),G_0,l)$. In consequence the gauge groupoid $\frac{P_0^u\times P_0^u}{U_0}\tto P_0^u/U_0$ is a subgroupoid of the gauge groupoid $\frac{P_0\times P_0}{U_0}\tto P_0/U_0$.

\begin{prop}\label{prop1U} The following isomorphism
 \unitlength=5mm \be\label{gaugethmU}\begin{picture}(11,4.6)
    \put(1,4){\makebox(0,0){$\frac{P^u_0\times P^u_0}{U_0}$}}
    \put(8,4){\makebox(0,0){$\U_{p_0}(\M)$}}
    \put(1,-1){\makebox(0,0){$P^u_0/U_0$}}
    \put(8,-1){\makebox(0,0){$\Ll_{p_0}(\M)$}}
    \put(1.2,3){\vector(0,-1){3}}
    \put(0.7,3){\vector(0,-1){3}}
    \put(8.2,3){\vector(0,-1){3}}
    \put(7.7,3){\vector(0,-1){3}}
    \put(3,4){\vector(1,0){3}}
    \put(2.7,-1){\vector(1,0){3.7}}
    \put(0.1,1.4){\makebox(0,0){$\Ss$}}
    \put(2.2,1.4){\makebox(0,0){$\Tt$}}
    \put(9.1,1.4){\makebox(0,0){$l$}}
    \put(6.8,1.4){\makebox(0,0){$r$}}
    \put(4.5,4.5){\makebox(0,0){$\phi^u$}}
    \put(4.5,-0.5){\makebox(0,0){$\varphi^u $}}
    \end{picture},\ee

    \bigskip
    of real Banach-Lie groupoids takes place, where
    \be \phi^u:\frac{P^u_0\times P^u_0}{U_0}\ \ni\ \langle\eta,\xi\rangle\mapsto \eta\xi^*\ \in\ \U_{p_0}(\M)\ee
    and
    \be \varphi^u:P^u_0/U_0\ \ni\ [\eta]\mapsto \eta\eta^*\ \in\ \Ll_{p_0}(\M).\ee
\end{prop}

From Proposition \ref{prop1U} we conclude that the groupoid of partial isometries $\U_{p_0}(\M)\tto\Ll_{p_0}(\M)$ is a subgroupoid of the groupoid $\G_{p_0}(\M)\tto\Ll_{p_0}(\M)$.

Let $(\vartheta,\eta)\in T_\eta P_0^u$, i.e. $\vartheta=\frac{d}{dt}\eta(t)|_{t=0}$, where values $\eta(t)=\gamma_t(\eta)$ of the flow  $\gamma_t(\eta)$  satisfy the condition (\ref{PU}). Thus one obtains
\be\label{polev}\eta^*\vartheta(\eta,\eta^*)+\vartheta(\eta,\eta^*)^*\eta=0,\ee
i.e. $\eta^*\vartheta(\eta,\eta^*)\in ip_0\M p_0$.

Now let us take $\mathfrak{X}\in \Gamma_{U_0}^\infty(TP_0^u)$. Similarly as in (\ref{Xeta}) the corresponding section $X\in \Gamma^\infty\A_{p_0}^u(\M)$ of the algebroid $\A_{p_0}^u(\M)$ of the groupoid $\U_{p_0}(\M)\tto\Ll_{p_0}(\M)$ is given by
\be\label{Xleta} X(l(\eta))=\vartheta(\eta,\eta^*)\eta^*.\ee

\begin{prop}\label{nawias4} If the sections $X_1, X_2 \in \Gamma^\infty\A_{p_0}(\M)$ satisfy the condition (\ref{Xleta}) on the unitary principal bundle $P_0^u\subset P_0$ then their algebroid bracket $[X_1,X_2]_\A$ also satisfies this condition on $P_0^u$.
\end{prop}

\prf{ Substituting $X=[X_1,X_2]_A$ into (\ref{polev}) we obtain
\be\label{dowod}\eta^*(X+X^*)\eta= \eta^*([X_2,X_1]+[X_1^*,X_2^*])\eta+\ee
$$+\eta^*\left(\langle\frac{\partial X_2}{\partial\eta}, X_1\eta\rangle-\langle\frac{\partial X_1}{\partial\eta}, X_2\eta\rangle+\langle\frac{\partial X_2}{\partial\eta^*},\eta^* X_1^*\rangle-
\langle\frac{\partial X_1}{\partial\eta^*},\eta^* X_2^*\rangle\right.$$
$$+\left.\langle\frac{\partial X_2^*}{\partial\eta^*}, \eta^* X_1^*\rangle-\langle\frac{\partial X_1^*}{\partial\eta^*}, \eta^*X_2^*\rangle+\langle\frac{\partial X_2^*}{\partial\eta}, X_1\eta\rangle-
\langle\frac{\partial X_1^*}{\partial\eta}, X_2\eta\rangle\right)\eta=$$
$$ =\eta^*\left([X_2,X_1]+[X_1^*,X_2^*]\right)\eta+$$
$$+\eta^*\left((X_1+X_1^*)X_2-(X_2+X_2^*)X_1+X_2^*(X_1+X_1^*)-X_1^*(X_1+X_1^*)\right)\eta+$$
$$+\langle\frac{\partial\eta^*(X_2+X_2^*)\eta}{\partial\eta},X_1\eta\rangle-
\langle\frac{\partial\eta^*(X_1+X_1^*)\eta}{\partial\eta},X_2\eta\rangle$$
$$+\langle\frac{\partial\eta^*(X_2+X_2^*)\eta}{\partial\eta^*},\eta^* X_1^*\rangle-\langle\frac{\partial\eta^*(X_1+X_1^*)\eta}{\partial\eta^*},\eta^* X_2^*\rangle=0.$$

The last equality in (\ref{dowod}) follows from $\eta^*(X_k+X_k^*)\eta=0$ for $k=1,2$.
}
From the above proposition one has
\begin{cor} The algebroid $\A_{p_0}^u(\M)$ is a transitive subalgebroid of $\A_{p_0}(\M)$. The anchor map $\mathbf{a}$ of $\A_{p_0}^u(\M)$ is the restriction of  $\ Tl:\Gamma^\infty(\A_{p_0}(\M))\to \Gamma^\infty(T\Ll_{p_0}(\M))$ to $\Gamma^\infty(\A_{p_0}^u(\M))$ and thus it is given by (\ref{anchorT}).
\end{cor}
\prf{ It follows from (\ref{polev}) that the algebroid $\A_{p_0}^u(\M)\to \Ll_{p_0}(\M)$ is a real Banach subbundle of the algebroid $\A_{p_0}(\M)\to \Ll_{p_0}(\M)$.  From the Proposition \ref{nawias4} we conclude  that if $X_1, X_2 \in \Gamma^\infty(\A_{p_0}^u(\M))$ then $[X_1, X_2] \in \Gamma^\infty(\A_{p_0}^u(\M))$.
}

Repeating the considerations similar to the ones for (\ref{atyah}) we obtain the Atiyah sequence
\be\label{atyahU} 0\to ip_0\M ^hp_0\times_{Ad_{U_0}} P^u_0\stackrel{\iota}{\hookrightarrow} TP^u_0/ {U_0}\stackrel{\mathbf{a}}{\rightarrow} T(P^u_0/U_0)\to 0\ee
of the unitary principal bundle $P_0^u(\Ll_{p_0}(\M),U_0,l)$, where $ip_0\M^hp_0$ is the real Banach space of the anti-hermitian elements of $p_0\M p_0$. So, due to the Proposition \ref{prop1U} we identify $TP^u_0/ {U_0}$ with the algebroid $\A^u_{p_0}(\M)$ of the groupoid $\U_{p_0}(\M)\tto\Ll_{p_0}(\M)$ as well as $T(P^u_0/U_0)$ with $T(\Ll_{p_0}(\M)$. The associated vector bundle $ip_0\M ^hp_0\times_{Ad_{U_0}} P^u_0$ is isomorphic with the vertical tangent bundle of $P_0^u(\Ll_{p_0}(\M),U_0,l)$. The section $\mathfrak{X}\in\Gamma^\infty(TP^u_0/ {U_0})\cong  \Gamma^\infty(\A^u_{p_0}(\M))$ is given by (\ref{X}) where $\vartheta$ satisfies the condition (\ref{polev}).

Substituting (\ref{eta}) and (\ref{wzor}) into (\ref{polev}) we obtain
\be\label{nazwa} (p+y_p)^*(a_p+(p+y_p)b_p)+(a_p+(p+y_p)b_p)^*(p+y_p)=0.\ee
So, if the operator valued functions $a_p$ and $b_p$, defining the vector field $\mathfrak{X}\in \Gamma^\infty_{G_0}(TP_0)$ from (\ref{frakX}), satisfy (\ref{nazwa}) on $P_0^u$ then $\mathfrak{X}$ is $U_0$-invariant and tangent to $P_0^u$.

\begin{prop}
The condition (\ref{polev}) is preserved by the bracket (\ref{bracket}). The same is valid for the condition (\ref{nazwa}) and the bracket (\ref{nawias}).
\end{prop}
\prf{ In  order to prove the first part of the proposition let us take $\vartheta_1$ and $\vartheta_2$ such that
$$\eta  ^*\vartheta_1+\vartheta_1^*\eta=0\quad {\rm and }\quad
\eta  ^*\vartheta_2+\vartheta_2^*\eta=0.$$
The vector field  $[\mathfrak{X}_1,\mathfrak{X}_2]$ is defined by (\ref{X}) where
$$\vartheta=\langle\frac{\partial \vartheta_2}{\partial \eta},\vartheta_1\rangle-\langle\frac{\partial \vartheta_1}{\partial \eta},\vartheta_2\rangle+ \langle\frac{\partial \vartheta_2}{\partial \eta^*},\vartheta_1^*\rangle-\langle\frac{\partial \vartheta_1}{\partial \eta^*},\vartheta_2^*\rangle.$$
Thus we have
$$\eta  ^*\vartheta+\vartheta^*\eta=\eta^*\left( \langle\frac{\partial \vartheta_2}{\partial \eta},\vartheta_1\rangle-\langle\frac{\partial \vartheta_1}{\partial \eta},\vartheta_2\rangle+ \langle\frac{\partial \vartheta_2}{\partial \eta^*},\vartheta_1^*\rangle-\langle\frac{\partial \vartheta_1}{\partial \eta^*},\vartheta_2^*\rangle\right)+$$
$$+\left( \langle\frac{\partial \vartheta_2^*}{\partial \eta},\vartheta_1\rangle-\langle\frac{\partial \vartheta_1^*}{\partial \eta},\vartheta_2\rangle+ \langle\frac{\partial \vartheta_2^*}{\partial \eta^*},\vartheta_1^*\rangle-\langle\frac{\partial \vartheta_1^*}{\partial \eta^*},\vartheta_2^*\rangle\right)\eta=$$

$$=\langle\frac{\partial }{\partial \eta}(\eta^*\vartheta_2+\vartheta_2^*\eta),\vartheta_1\rangle-\langle\vartheta_2^*,\vartheta_1\rangle-\langle\frac{\partial }{\partial \eta}(\eta^*\vartheta_1+\vartheta_1^*\eta),\vartheta_2\rangle+\langle\vartheta_1^*,\vartheta_2\rangle+$$

$$+\langle\frac{\partial }{\partial \eta^*}(\eta^*\vartheta_2+\vartheta_2^*\eta),\vartheta_1^*\rangle-\langle\vartheta_1^*,\vartheta_2\rangle-\langle\frac{\partial }{\partial \eta^*}(\eta^*\vartheta_1+\vartheta_1^*\eta),\vartheta_2^*\rangle + \langle\vartheta_2^*,\vartheta_1\rangle=0.$$

We have the following equalities 
$$ \langle \frac{\partial b_{pi}}{\partial y_p},a_{pj}\rangle=-(p+y_p)^{-1}a_j((p+y_p)^*(p+y_p))^{-1}Y_{pi}+$$
\be\label{d1} +((p+y_p)^*(p+y_p))^{-1}\langle\frac{\partial Y_{pi}}{\partial y_p},a_{pj}\rangle+\ee
$$+(p+y_p)^{-1}a_j(p+y_p)^{-1}a_i-(p+y_p)^{-1}\langle\frac{\partial a_{pj}}{\partial y_p},a_{pi}\rangle $$
and
\be\label{d2} \langle \frac{\partial b_{pi}}{\partial y^*_p},a_{j}\rangle=-((p+y_p)^*(p+y_p))^{-1}a^*_{pj}(p+y_p)^{-1}Y_{pi}+
\ee
 $$+((p+y_p)^*(p+y_p))^{-1}\langle\frac{\partial Y_{pj}}{\partial y^*_p},a^*_{pj}\rangle-(p+y_p)^{-1}\langle\frac{\partial a_{pj}}{\partial y^*_{pi}},a^*_{pj}\rangle,$$
where  \be\label{ypi}Y_{pi}:=(p+y_{p})^*(a_{pi}+(p+y_p)b_{pi}).\ee

Substituting  (\ref{d1}), (\ref{d2}), into (\ref{ba})  we get
\be\label{bbp} b_p= \langle\frac{\partial b_{2p}}{\partial y_p},a_{1p}\rangle-\langle\frac{\partial b_{1p}}{\partial y_p},a_{2p}\rangle+\langle\frac{\partial b_{2p}}{\partial y^*_p},a^*_{1p}\rangle-\langle\frac{\partial b_{1p}}{\partial y^*_p},a^*_{2p}\rangle+[b_{1p},b_{2p}]=\ee
$$=((p+y_p)^*(p+y_p))^{-1}\left[ Y_{2p}((p+y_p)^*(p+y_p))^{-1}Y_{1p}-Y_{1p}((p+y_p)^*(p+y_p))^{-1}Y_{2p}-\right.$$
$$-Y_{2p}(p+y_p)^{-1}a_{1p}-
a_{1p}^*((p+y_p)^*)^{-1}Y_{2p}+a_{2p}^*((p+y_p)^*)^{-1}Y_{1p}+Y_{1p}(p+y_p)^{-1}a_{2p}+$$
$$\left.
+\langle\frac{\partial Y_{2p}}{\partial y_p},a_{1p}\rangle+\langle\frac{\partial Y_{2p}}{\partial y_p^*},a_{1p}^*\rangle-
\langle\frac{\partial Y_{1p}}{\partial y_p},a_{1p}\rangle-\langle\frac{\partial Y_{1p}}{\partial y_p^*},a_{1p}^*\rangle\right]$$
$$-(p+y_p)^{-1}\left[\langle\frac{\partial a_{2p}}{\partial y_p},a_{1p}\rangle-\langle\frac{\partial a_{1p}}{\partial y_p},a_{2p}\rangle+\langle\frac{\partial a_{2p}}{\partial y^*_p},a^*_{1p}\rangle-\langle\frac{\partial a_{1p}}{\partial y^*_p},a^*_{2p}\rangle\right].$$
From (\ref{bbp}) we obtain
$$ Y_{p}:=(p+y_{p})^*(a_{p}+(p+y_p)b_{p})=Y_{2p}((p+y_p)^*(p+y_p))^{-1}Y_{1p}-$$
$$-Y_{1p}((p+y_p)^*(p+y_p))^{-1}Y_{2p}-Y_{2p}(p+y_p)^{-1}a_{1p}-$$
$$-a_{1p}^*((p+y_p)^*)^{-1}Y_{2p}+a_{2p}^*((p+y_p)^*)^{-1}Y_{1p}+Y_{1p}(p+y_p)^{-1}a_{2p}+$$
\be\label{yyp} +\langle\frac{\partial Y_{2p}}{\partial y_p},a_{1p}\rangle+\langle\frac{\partial Y_{2p}}{\partial y_p^*},a_{1p}^*\rangle-
\langle\frac{\partial Y_{1p}}{\partial y_p},a_{1p}\rangle-\langle\frac{\partial Y_{1p}}{\partial y_p^*},a_{1p}^*\rangle.\ee
Since  $Y_{pi}+Y_{pi}^*=0$, for $i=1,2$, it follows from (\ref{yyp})  that  $Y_{p}+Y_{p}^*=0$, i.e. the condition (\ref{nazwa}) is fulfilled for $[\mathfrak{X}_1,\mathfrak{X}_2 ]$.
}

\section{Example}
In this section we discuss the case  $\M=L^\infty(\H)$, where $\H$ is a separable complex Hilbert space with a fixed orthonormal basis $\{|e_i\rangle\}_{i=1}^\infty$. Using Dirac notation we set \be\label{elem1} p_0=\proj{e_1}{e_1}+\proj{e_2}{e_2}+...+\proj{e_{N}}{e_{N}}.\ee
Hence for $\vartheta\in (L^\infty(\H))p_0$ and for $\eta\in P_0\subset L^\infty(\H)$ we have
\be\label{elem2} \vartheta=\proj{\vartheta_1}{e_1}+\proj{\vartheta_2}{e_2}+...+\proj{\vartheta_{N}}{e_{N}},\ee
where $\vartheta_1,...,\vartheta_N\in \H$, and
\be\label{elem3} \eta=\proj{\eta_1}{e_1}+\proj{\eta_2}{e_2}+...+\proj{\eta_{N}}{e_{N}},\ee
where $\eta_1,...,\eta_N\in \H$ are linearly independent vectors. The polar decomposition
$$ \eta=u|\eta|$$ of $\eta\in P_0$ consists
\be\label{elem4} |\eta|=(\eta^*\eta)^{\frac{1}{2}}=\left(\sum _{k,l=1}^N\ket{e_k}\sc{\eta_k}{\eta_l}\bra{e_l}\right)^{\frac{1}{2}}\ee
and
\be u=\eta(\eta^*\eta)^{-\frac{1}{2}}=\proj{u_1}{e_1}+\proj{u_2}{e_2}+...+\proj{u_{N}}{e_{N}},\ee
where $\sc{u_k}{u_l}=\delta_{kl}.$

Any $x\in \G_{p_{0}}(L^\infty(\H))$ one can express in the following way
\be x=\proj{\eta_1}{\xi_1}+...+\proj{\eta_N}{\xi_N}\ee
where $\eta_1,...,\eta_N\in\H$ and $\xi_1,...,\xi_N\in\H$ are sequences of linearly independent vectors. Hence
\be r(x)=\proj{u_1}{u_1}+\proj{u_2}{u_2}+...+\proj{u_{N}}{u_{N}}\ee
and\
\be l(x)=\proj{v_1}{v_1}+...+\proj{v_{N}}{v_{N}},\ee
where
\be v=\xi(\xi ^*\xi)^{-1}=\proj{v_1}{e_1}+...+\proj{v_{N}}{e_{N}}\ee
and $\xi=v|\xi|$.

The groupoid isomorphism described in Proposition \ref{prop22} is given by
\be \phi(\langle\eta,\xi\rangle)=\eta\xi^{-1}=\sum_{k,l=1}^N\ket{\eta_k}\sc{e_k}{(\eta^*\eta)^{-1}e_l}\bra{\xi_l}\ee
and
\be \varphi([\eta])=\eta\eta^{-1}=\proj{u_1}{u_1}+...+\proj{u_{N}}{u_{N}}.\ee
Since of (\ref{elem3}) we will identify $\eta$ and $\eta^*$ with $(\eta_1,...,\eta_N)\in \H\times\cdots\times\H$ and $(\ol\eta_1,...,\ol\eta_N)\in \ol\H\times\cdots\times\ol\H$, respectively. By $\ol\H$ we denoted the Hilbert space conjugated to $\H$. Using above convention we express the algebroid bracket (\ref{bracket})  as follows
\be [\mathfrak{X}_1,\mathfrak{X}_2]=\ee
$$=\sum _{n,k=1}^N\left\langle\frac{\partial}{\partial \eta_n},\langle\frac{\partial\vartheta_{2n}}{\partial\eta_k},\vartheta_{1k}\rangle-\langle\frac{\partial\vartheta_{1n}}{\partial\eta_k},\vartheta_{2k}\rangle+
\langle\frac{\partial\vartheta_{2n}}{\partial\ol\eta_k},\ol\vartheta_{1k}\rangle-
\langle\frac{\partial\vartheta_{1n}}{\partial\ol\eta_k},\ol\vartheta_{2k}\rangle\right\rangle+$$
$$+\sum _{n,k=1}^N\left\langle\frac{\partial}{\partial \ol\eta_n},\langle\frac{\partial\ol\vartheta_{2n}}{\partial\eta_k},\vartheta_{1k}\rangle-\langle\frac{\partial\ol\vartheta_{1n}}{\partial\eta_k},\vartheta_{2k}\rangle+
\langle\frac{\partial\ol\vartheta_{2n}}{\partial\ol\eta_k},\ol\vartheta_{1k}\rangle-
\langle\frac{\partial\ol\vartheta_{1n}}{\partial\ol\eta_k},\ol\vartheta_{2k}\rangle\right\rangle,$$
where $\mathfrak{X}_l$ is given by
\be \mathfrak{X}_l=\sum _{n=1}^N\left(\langle\frac{\partial}{\partial \eta_n},\vartheta_{ln}\rangle+\langle\frac{\partial}{\partial \ol\eta_n},\ol\vartheta_{ln}\rangle\right)\quad {\rm for}\quad l=1,2.\ee

Now let us present $(y_p,z_p)$-coordinates description of the groupoid $\G_{p_{0}}(L^\infty(\H))\tto \Ll_{p_{0}}(L^\infty(\H))$. For this reason let us consider projections given by
\be p_{n_1n_2...n_N}:=\proj{e_{n_1}}{e_{n_1}}+\proj{e_{n_2}}{e_{n_2}}+...+\proj{e_{n_N}}{e_{n_N}},\ee
where $n_1<n_2<\cdots <n_N$. We observe that the domains $\Pi_{n_1n_2...n_N}:=\Pi_{p_{n_1n_2...n_N}}$ and $l^{-1}(\Pi_{n_1n_2...n_N})$ of the charts
$$\varphi_{n_1n_2...n_N}:\Pi_{n_1n_2...n_N}\to (1-p_{n_1n_2...n_N})(L^\infty(\H))p_{n_1n_2...n_N}$$
and
$$\psi_{n_1n_2...n_N}:l^{-1}(\Pi_{n_1n_2...n_N})\to(1-p_{n_1n_2...n_N})(L^\infty(\H))p_{n_1n_2...n_N}\oplus p_{n_1n_2...n_N}(L^\infty(\H))p_0$$
cover $\Ll_{p_0}(L^\infty(\H))$ and  $P_0$, respectively.
Using (\ref{yetazeta}) we find that
\be z_ {n_1n_2...n_N}=p_{n_1n_2...n_N}\eta=\sum_{r=1}^N
\ket{e_{n_r}}\sc{e_{n_r}}{\eta_s}\bra{e_s}\ee
and
\be y_{n_1n_2...n_N}=(1-p_{n_1n_2...n_N})\eta(p_{n_1n_2...n_N}\eta)^{-1}=\ee
$$=\left(\sum_{ r\not=n_1,.,n_N}^\infty\proj{e_r}{e_r}\right)\left(\sum_{k=1}^N\proj{\eta_k}{e_k}\right)
\left(\sum_{s,t=1}^N\proj{e_s}{e_s}z^{-1}_{n_1n_2...n_N}\proj{e_{n_t}}{e_{n_t}}\right)=$$
$$=\sum_{r\not=n_1,.,n_N}^\infty\sum_{t=1}^N\ket{e_r}\left(\sum_{k=1}^N\sc{e_r}{\eta_k}
\bra{e_k}z^{-1}_{n_1n_2...n_N}\ket{e_{n_t}}\right)\bra{e_{n_t}}=$$
$$=\sum_{r\not=n_1,.,n_N}^\infty\sum_{t=1}^N\ket{e_r}\left(\sum_{k=1}^N\sc{e_r}{\eta_k}
\bra{e_k}z^{-1}_{n_1n_2...n_N}\ket{e_{n_t}}\right)\bra{e_{n_t}},$$
where
$$z^{-1}_{n_1n_2...n_N}:p_{n_1n_2...n_N}\H\to p_0\H$$ is the inverse of the isomorphism $z_{n_1n_2...n_N}:p_0\H\to p_{n_1n_2...n_N}\H$ of the $N$-dimensional Hilbert subspaces of $\H$. The matrix element of the operators $z_{n_1n_2...n_N}$ and $y_{n_1n_2...n_N}$ are the following:
\be z_{n_1n_2...n_N}^{ts}=\sc{e_t}{\eta_s}\ee
where $t\in\{n_1,n_2,...,n_N\}, \ s\in \{1,2,...,N\}$ and
\be y_{n_1n_2...n_N}^{rt}=\sum_{k=1}^N\sc{e_r}{\eta_k}\bra{e_k}z^{-1}_{n_1n_2...n_N}\ket{e_{t}}\ee
where $r\in \mathbb{N}\setminus \{n_1,n_2,...,n_N\}$. In the same way we define
\be b_{n_1n_2...n_N}^{ts}:=\bra{e_{t}}b_{n_1n_2...n_N}\ket{e_{s}}\ee
and
\be a_{n_1n_2...n_N}^{rt}:=\bra{e_{r}}a_{n_1n_2...n_N}\ket{e_{t}}\ee
the matrix elements of $b_{n_1n_2...n_N}$ and $a_{n_1n_2...n_N}$. Now, passing in (\ref{nawias}), (\ref{ap}) and (\ref{bp}) to matrix element description we obtain
\be [\mathfrak{X}_1,\mathfrak{X}_2]=\sum_{r\not=n_1,.,n_N}^\infty\sum_{t=n_1,.,n_N}\left(a_{n_1...n_N}^{rt} \frac{\partial }{\partial y^{rt}_{n_1...n_N}}+\ol a_{n_1...n_N}^{rt} \frac{\partial }{\partial \ol y^{rt}_{n_1...n_N}}\right)+\ee
$$+\sum_{m=n_1,.,n_N}\sum_{s,k=1,.,N}\left( b_{n_1...n_N}^{ms}z^{n_sk}_{n_1...n_N} \frac{\partial }{\partial z^{mk}_{n_1...n_N}}+\ol z^{n_sk}_{n_1...n_N}\ol b_{n_1...n_N}^{ms} \frac{\partial }{\partial \ol z^{mk}_{n_1...n_N}}\right),$$
where
$$ b^{ts}_{n_1...n_N}=\sum_{j\not=n_1,.,n_N}^\infty\sum_{m=n_1,.,n_N}\left(a_{1{n_1...n_N}}^{jm}\frac{\partial b_{2{n_1...n_N}}^{ts}}{\partial y^{jm}_{n_1...n_N}}-
a_{2{n_1...n_N}}^{jm}\frac{\partial b_{1{n_1...n_N}}^{ts}}{\partial y^{jm}_{n_1...n_N}}\right)+$$
$$+\sum_{j\not=n_1,.,n_N}^\infty\sum_{m=n_1,.,n_N}\left(
\ol a_{1{n_1...n_N}}^{jm}\frac{\partial b_{2{n_1...n_N}}^{ts}}{\partial \ol y^{jm}_{n_1...n_N}}-
\ol a_{2{n_1...n_N}}^{jm}\frac{\partial b_{1{n_1...n_N}}^{ts}}{\partial \ol y^{jm}_{n_1...n_N}}\right)+$$
\be+\sum_{m=1,.,N
}\left( b_{2{n_1...n_N}}^{tm}b_{1{n_1...n_N}}^{n_ms}-b_{1{n_1...n_N}}^{tm}b_{2{n_1...n_N}}^{n_ms}\right)\ee

and

\be a^{rt}_{n_1...n_N}=\sum_{j\not=n_1,.,n_N}^\infty\sum_{r=n_1,.,n_N}\left(a_{1{n_1...n_N}}^{jm}\frac{\partial a_{2{n_1...n_N}}^{rt}}{\partial y^{jm}_{n_1...n_N}}-
a_{2{n_1...n_N}}^{jm}\frac{\partial a_{1,{n_1...n_N}}^{rt}}{\partial y^{jm}_{n_1...n_N}}\right)+\ee
$$+\sum_{j\not=n_1,.,n_N}^\infty\sum_{r=n_1,.,n_N}\left(\ol a_{1{n_1...n_N}}^{jm}\frac{\partial a_{2{n_1...n_N}}^{rt}}{\partial \ol y^{jm}_{n_1...n_N}}-
\ol a_{2{n_1...n_N}}^{jm}\frac{\partial a_{1{n_1...n_N}}^{rt}}{\partial \ol y^{jm}_{n_1...n_N}}\right).$$

The anchor map  in coordinates $(y_{n_1n_2...n_N},
\ol y_{n_1n_2...n_N})$ takes the following form

\be \mathbf{a}(\mathfrak{X})=\sum_{r\not=n_1,.,n_N}^\infty\sum_{j=n_1,.,n_N}
 \left( a_{n_1...n_N}^{rj} \frac{\partial }{\partial y^{rj}_{n_1...n_N}}+\ol a_{n_1...n_N}
^{rj}
 \frac{\partial }{\partial \ol y^{rj}_{n_1...n_N}}\right)\ee
for $\mathfrak{X}\in \Gamma^\infty \A_{p_0}(\M)$.

Ending the example let us mention that $P_0$ is the Stiefel principal bundle of $N$-frames of the tautological vector bundle $\mathbb{E}\to G(N,\H)$ over the Grassmannian $G(N,\H)$ of $N$-dimensional subspaces of the Hilbert space $\H$. On the other hand groupoids $\G_{p_{0}}(L^\infty(\H))\tto \Ll_{p_{0}}(L^\infty(\H))\cong G(N,\H)$ and $\U_{p_{0}}(L^\infty(\H))\tto \Ll_{p_{0}}(L^\infty(\H))\cong G(N,\H)$ are frame groupoid and unitary frame groupoid of $\mathbb{E}\to G(N,\H)$, respectively.

\section{Other related groupoids}

 Following \cite{brown} we present a construction of a quotient groupoid which will be useful in the subsequent considerations. Let $H\tto B$ be a wide subgroupoid  of a groupoid $G\tto B$ such that $H_a^b$, where $a,b\in B$, has at most one element. By $B/ H$ we denote the set of orbits of the action $H$ on the base manifold $B$, i.e. $a,b\in B$ belong to the same orbit $[a]$ if $h\in H$ exists such that $\Ss(h)=a$ and $\Tt(h)=b$. Note that if such $h\in H$ exists it is unique.

 We define the quotient set $G/ H$ in the following way. Elements  $g_1,g_2\in G$ are equivalent if and only if there exist $h_1,h_2\in H$ such that $g_2=h_2 g_1 h_1$. The product of $[g], [f]\in G/ H$, where $\Ss([g])=\Tt([f])$ we define as follows
 \be [g]\cdot[f]:=[g h f],\ee
 where $h\in H_{\Tt(f)}^{\Ss(g)}$
 and other structure maps are defined by:
 $$\ \Ss([g]):=[\Ss(g)],\  \ \ \Tt([g]):=[\Tt(g)],\ \  \ [g]^{-1}:=[g^{-1}],\ \  \ \varepsilon([g]):=[\varepsilon(g)].$$

 \begin{prop}\label{qoutient} Let $H\tto B$ be  such a  wide subgroupoid of the groupoid $G\tto B$ that for any $a,b\in B$ the set $H_a^b$ has at most one element. Then on the quotient space $G/ H$ there exists the groupoid structure $G/ H\tto B/ H$.

 \end{prop}

 \begin{example}
 Let us consider a principal bundle $P(M,G,\pi)$. We take the pair groupoid $P\times P\tto P$ and action groupoid $G\sphericalangle P\tto P$ is a subgroupoid of $P\times P\tto P$. One easily verifies that in this case the assumptions of the proposition are fulfilled and the quotient groupoid $\frac{P\times P}{G\sphericalangle P}\tto P/{G}$ is isomorphic to the pair groupoid $B\times B\tto B$.
 \end{example}

 The construction of the quotient groupoid $G/ H\tto B/ H$ is in a sense complementary to the construction of the quotient groupoid $G/ N\tto B$, where $N\tto B$ is the normal subgroupoid of $G\tto B$, see e.g. \cite{mac}. Thereby one does not change the base set but replaces arrows by their equivalence classes. In our case since the action of $H\tto B$ on the base $B$ is free, the equivalence relation on $B$ determines equivalence relation on the groupoid set $G$.

 By structural (frame) groupoid $\G^{lin}\mathbb{E}\tto B$ of a Banach vector bundle $\pi:E\to B$ we will mean the groupoid  which consists of continuous linear fibre isomorphisms $G_{ba}:\mathbb{E}_a\to\mathbb{E}_b$. All structural maps of $\G^{lin}\mathbb{E}\tto B$ are defined in the same way as in the finite dimensional case, e.g. the product of $G_{cb}$ and  $G_{ba}$ is given by the superposition $G_{cb} G_{ba}:\mathbb{E}_a\to\mathbb{E}_c$ of isomorphisms (e.g. see \cite{mac}).

 Among other groupoids we will consider the following frame groupoids $\G^{lin}\A_{p_0}(\M)\tto\Ll_{p_0}(\M)$, $\G^{lin}\M^R_{p_0}(\M)\tto\Ll_{p_0}(\M)$ and $\G^{lin}TP_0\tto P_0$. We note firstly that the groupoid $\G^{lin}TP_0\tto P_0$ is isomorphic with the trivial groupoid $P_0\times Aut(\M p_0)\times P_0\tto P_0$. Recall that the product of $(\lambda, \Lambda, \xi), (\xi,\Gamma,\eta)\in P_0\times Aut(\M p_0)\times P_0$ is given by
 $$(\lambda, \Lambda, \xi)(\xi,\Gamma,\eta)=(\lambda, \Lambda\Gamma,\eta).$$ The other groupoid maps are defined by
 $$(\lambda, \Lambda, \xi)^{-1}=(\xi, \Lambda^{-1}, \lambda),$$
 $$\Ss(\lambda, \Lambda, \xi)=\xi$$
 $$\Tt(\lambda, \Lambda, \xi)=\lambda$$
 and $\varepsilon:P_0\to P_0\times Aut(\M p_0)\times P_0$ by
 $$\varepsilon(\eta)=(\eta,id,\eta).$$
 Since the tangent vector bundle $TP_0\to P_0$ is trivial and isomorphic  to the bundle $\M p_0\times P_0\to P_0$ one defines the groupoid isomorphism
 \be I:P_0\times Aut(\M p_0)\times P_0\ \ni \ (\xi,\Gamma,\eta)\mapsto \Gamma_{\xi\eta}\  \in\ \G^{lin}TP_0\ee
between $\G^{lin}TP_0\tto P_0$ and $P_0\times Aut(\M p_0)\times P_0\tto P_0$ by
  \be\Gamma_{\xi\eta}(\vartheta,\eta):=(\Gamma(\vartheta),\xi),\ee
 where $(\vartheta,\eta)\in T_\eta P_0=\M p_0\times\{\eta\}$ and $(\Gamma(\vartheta),\xi)\in T_\xi P_0=\M p_0\times\{\xi\}$.

 We note that action groupoid $G_0\sphericalangle P_0\tto P_0$ can be considered as a subgroupoid of $P_0\times Aut(\M p_0)\times P_0\tto P_0$, i.e.
 \be G_0\sphericalangle P_0\ \ni\ (g,\eta)\to (\eta g, R_g,\eta)\ \in P_0\times Aut(\M p_0)\times P_0,\ee
 where $R_g\in Aut(\M p_0)$ is defined by \be R_g \vartheta:=\vartheta g.\ee

 We also consider the subgroupoid $P_0\times Aut^{mod}(\M p_0)\times P_0\tto P_0$ of the trivial groupoid  consisting of such elements $(\xi,\Gamma,\eta)\in P_0\times Aut(\M p_0)\times P_0$ that
 \be\label{Gamma} \Gamma(x\vartheta)=x\Gamma\vartheta\ee
 for $x\in \M$ and $\vartheta\in \M p_0$. We will denote this subgroupoid by $P_0\times Aut^{mod}(\M p_0)\times P_0\tto P_0$, where by $P_0\times Aut^{mod}(\M p_0)\times P_0\tto P_0$ we denote the group of automorphisms of the $\M$-right modules $\M p_0$. The action groupoid $G_0\sphericalangle P_0\tto P_0$ is also a subgroupoid of $P_0\times Aut^{mod}(\M p_0)\times P_0\tto P_0$. It can by easily seen that the conditions of the Proposition \ref{qoutient} are fulfilled for $G_0\sphericalangle P_0\tto P_0$   considered as a subgroupoid of groupoids $P_0\times Aut(\M p_0)\times P_0\tto P_0$ and $P_0\times Aut^{mod}(\M p_0)\times P_0\tto P_0$. Therefore we have the following  proposition.

 \begin{prop}
 One has the groupoid isomorphisms:

 \unitlength=5mm \be\label{iso1}\begin{picture}(11,4.6)
    \put(-1,4){\makebox(0,0){$\frac{P_0\times Aut(\M p_0)\times P_0}{G_0\sphericalangle P_0}$}}
    \put(9,4){\makebox(0,0){$\G^{lin}\A_{p_0}(\M)$}}
    \put(-1,-1){\makebox(0,0){$\Ll_{p_0}(\M)$}}
    \put(8,-1){\makebox(0,0){$\Ll_{p_0}(\M)$}}
    \put(-1.2,3){\vector(0,-1){3}}
    \put(-0.7,3){\vector(0,-1){3}}
    \put(8.2,3){\vector(0,-1){3}}
    \put(7.7,3){\vector(0,-1){3}}
    \put(2.5,4){\vector(1,0){3.7}}
    \put(0.7,-1){\vector(1,0){5.7}}
    \put(4.5,4.5){\makebox(0,0){$[\phi]$}}
    \put(3.5,-0.5){\makebox(0,0){$id $}}
    \end{picture},\ee
    \newline
    \bigskip
 and
 \unitlength=5mm \be\label{iso2}\begin{picture}(11,4.6)
    \put(-1,4){\makebox(0,0){$\frac{P_0\times Aut^{mod}(\M p_0)\times P_0}{G_0\sphericalangle P_0}$}}
    \put(8,4){\makebox(0,0){$\G_{p_0}(\M)$}}
    \put(-1,-1){\makebox(0,0){$\Ll_{p_0}(\M)$}}
    \put(8,-1){\makebox(0,0){$\Ll_{p_0}(\M)$}}
    \put(-1.2,3){\vector(0,-1){3}}
    \put(-0.7,3){\vector(0,-1){3}}
    \put(8.2,3){\vector(0,-1){3}}
    \put(7.7,3){\vector(0,-1){3}}
    \put(2.9,4){\vector(1,0){3.3}}
    \put(0.7,-1){\vector(1,0){5.7}}
    \put(4.5,4.5){\makebox(0,0){$[\phi]$}}
    \put(3.5,-0.5){\makebox(0,0){$id $}}
    \end{picture},\ee
    \newline

  \end{prop}

  \prf{Let us define the groupoid map
  \unitlength=5mm \be\label{iso3}\begin{picture}(11,4.6)
    \put(-1,4){\makebox(0,0){$P_0\times Aut(\M p_0)\times P_0$}}
    \put(9,4){\makebox(0,0){$\G^{lin}\A_{p_0}(\M)$}}
    \put(-1,-1){\makebox(0,0){$P_0$}}
    \put(8,-1){\makebox(0,0){$\Ll_{p_0}(\M)$}}
    \put(-1.2,3){\vector(0,-1){3}}
    \put(-0.7,3){\vector(0,-1){3}}
    \put(8.2,3){\vector(0,-1){3}}
    \put(7.7,3){\vector(0,-1){3}}
    \put(2.9,4){\vector(1,0){3.3}}
    \put(0.7,-1){\vector(1,0){5.7}}
    \put(4.5,4.4){\makebox(0,0){$\phi$}}
    \put(3.5,-0.5){\makebox(0,0){$l $}}
    \end{picture},\ee
    \newline
    where $\phi_{(\xi,\Gamma,\eta)}$ acts on $[(\vartheta,\eta)]\in \A_{p_0}(\M)$ in the following way
    \be \phi_{(\xi,\Gamma,\eta)}[(\vartheta,\eta)]:=[(\Gamma(\vartheta),\xi)].\ee We note that $[(\vartheta,\eta)]\in\pi^{-1}(l(\eta))$ and $[(\Gamma(\vartheta),\xi)]\in\pi^{-1}(l(\xi))$. Any bounded linear isomorphism $\Gamma_{l(\xi)l(\eta)}:\pi^{-1}(l(\eta))\to\pi^{-1}(l(\xi))$ between the fibres of $\pi:\A_{p_0}\to\Ll_{p_0}(\M)$ is given by
    \be \Gamma_{l(\xi)l(\eta)}([\vartheta,\eta)])=[(\Gamma\vartheta,\xi)],\ee
    where $\Gamma\in Aut(\M p_0)$ is determined up to the transformation
    \be \Gamma '=R_h\Gamma R_{g^{-1}},\quad g,h\in G_0.\ee
    So, $\phi$ maps $P_0\times Aut(\M p_0)\times P_0$ onto $\G^{lin}\A_{p_0}(\M)$. The commutativity of $\phi$  and $l$ with sources and target maps in (\ref{iso3}) is obvious. The groupoid product
    \be \phi_{(\lambda,\Lambda,\xi)}\circ\phi_{(\xi,\Gamma,\eta)}=\phi_{(\lambda,\Lambda\Gamma,\eta)}\ee is also preserved by $\phi$.
    For any $(\eta g,R_g,\eta), \ (\xi h,R_h,\xi)\in P_0\times Aut(\M p_0)\times P_0$ we have
    \be \phi_{(\xi h,R_h \Gamma R_{g^{-1}},\eta g)}=\phi_{(\xi h,R_h ,\xi)}\circ \phi_{(\xi, \Gamma,\eta g)}\circ \phi_{(\eta g,R_g,\eta)^{-1}}.\ee
    From the above facts it follows that (\ref{iso3}) defines a morphism of groupoids. After the factorization by subgroupoid $G_0\sphericalangle P_0\tto P_0$ we obtain (\ref{iso1}).

    From (\ref{Gamma}) it follows that $(\xi, \Gamma,\eta)\in P_0\times Aut^{mod}(\M p_0)\times P_0$ if and only if $\Gamma=R_{\Gamma(\vartheta)}$. So, one can identify $P_0\times Aut^{mod}(\M p_0)\times P_0\tto P_0$ with the trivial groupoid $P_0\times G_0\times P_0\tto P_0$. For $(\xi h,h,\xi),(\eta g,g,\eta)\in G_0\sphericalangle P_0$ and $(\xi, f, \eta)\in P_0\times G_0\times P_0$ one has
    \be (\xi h,h,\xi)(\xi, f, \eta)(\eta g,g,\eta)^{-1}=(\xi h,hfg^{-1},\eta g).\ee
    The above enables one to obtain the following sequence of equalities
    $$ [(\xi,f,\eta)]=[(\xi h,hfg^{-1},\eta g)]=[(\xi g,gfg^{-1},\eta g)]=[(\xi f,f,\eta f)]\cong $$
    \be\label{warstwy}\cong [(\xi f,\eta f)]=[(\xi,\eta)]\ee
    for $[(\xi,f,\eta)]\in \frac{P_0\times G_0\times P_0}{G_0\sphericalangle P_0}$ and $[(\xi,\eta)]\in \frac{P_0\times P_0}{G_0}\cong \G_{p_0}(\M)$. From (\ref{warstwy}) we derive (\ref{iso2}).}

\section{Algebroid of the groupoid $\G^{lin}\A_{p_0}(\M)\tto\Ll_{p_0}(\M)$}

In what follows we will denote the algebroid of $\G^{lin}\A_{p_0}(\M)\tto\Ll_{p_0}(\M)$ by $\A^{lin}_{p_0}(\M)$ and describe its structure  in terms of the trivial groupoid $P_0\times Aut(\M p_0)\times P_0\tto P_0$.

Let us start noting that a smooth bisection $\sigma$ of the trivial groupoid is defined as follows:
\be\label{sigma} \sigma(\eta):=(\gamma(\eta),\Gamma(\eta),\eta)\ee
where $\gamma:P_0\to P_0$ is a diffeomorphism of $P_0$ and $\Gamma:P_0\to Aut(\M p_0)$ is a smooth map of $P_0$ in the group of automorphisms of the complex Banach space $\M p_0$. The product $\sigma_1\ast\sigma_2$ of two bisections of the trivial groupoid $P_0\times Aut(\M p_0)\times P_0\tto P_0$ is given by
\be\label{sigma2}(\sigma_1\ast\sigma_2)(\eta):=\left((\gamma_1\circ \gamma_2)(\eta), \Gamma_1(\gamma_2(\eta))\Gamma_2(\eta),\eta\right).\ee
Any bisection $\sigma$ defines linear continuous isomorphism $\Phi_\sigma(
\eta):T_\eta P_0\to T_{\gamma(\eta)} P_0$ of the fibres
\be\label{sigma3} \Phi_\sigma(\vartheta,\eta):=(\Gamma(\eta)\vartheta,\gamma(\eta)),\ee
of the tangent bundle $TP_0$, where $(\vartheta,\eta)\in T_\eta P_0$.

We will consider a one-parameter subgroup of bisections $\sigma_{t+s}=\sigma_t\ast\sigma_s$, i.e.
\be\label{bis1}\gamma_{t+s}=\gamma_t\circ\gamma_s\ee
\be\label{bis2} \Gamma_{t+s}(\eta)=\Gamma_t(\gamma_s(\eta))\Gamma_s(\eta),\ee
such that maps $\vartheta:P_0\to TP_0$ and $\Theta:P_0\to End(\M p_0)$, defined by
\be\label{pole1} \vartheta(\eta,\eta^*):=\frac{d}{ds}\gamma_s(\eta)|_{s=0}=\frac{d}{ds}\gamma(s,\eta,\eta^*)|_{s=0},\ee
\be\label{pole2} \Theta(\eta,\eta^*):=\frac{d}{ds}\Gamma_s(\eta)|_{s=0}=\frac{d}{ds}\Gamma(s,\eta,\eta^*)|_{s=0},\ee
are a smooth maps. Let us note that $\gamma:P_0\to P_0$  and $\Gamma:P_0\to Aut(\M p_0)$ as well as $\gamma_t:P_0\to P_0$  and $\Gamma_t:P_0\to Aut(\M p_0)$ are not holomorphic maps in general. So, in (\ref{pole1}),  (\ref{pole2}) and in formulae presented below we use the following notation:
\be \gamma(\eta)=\gamma(\eta,\eta^*) \quad {\rm and} \quad \Gamma(\eta)=\Gamma(\eta,\eta^*).\ee
\be \gamma_s(\eta)=\gamma(s,\eta,\eta^*) \quad {\rm and} \quad \Gamma_s(\eta)=\Gamma(s,\eta,\eta^*).\ee

Let a section $\mathcal{Z}\in \Gamma^\infty(TP_0)$ be given by  $ \mathcal{Z}(\eta)=(w(\eta,\eta^*),\eta)$. Then one defines the flow
$\Sigma_t: \Gamma^\infty(TP_0)\to \Gamma^\infty(TP_0)$, \ $t\in \mathbb{R}$, in the following way
\be \label{flow} (\Sigma_t\mathcal{Z})(\eta):=\Phi_{\sigma_t}(\mathcal{Z}(\gamma_{-t}(\eta)))=
\left(\Gamma_t(\gamma_{-t}(\eta))(w(\gamma_{-t}(\eta),\gamma_{-t}(\eta)^*)),\eta\right) .
\ee

Differentiating both sides of (\ref{flow}) with respect to the parameter $t\in \mathbb{R}$ we obtain a first-order differential operator

\be\label{diffop}(D\mathcal{Z})(\eta):=\left(((-\mathfrak{X}+\Theta)w)(\eta,\eta^*),\eta\right)\ee
which is a derivation of the module of sections $\Gamma^\infty(TP_0)$, i.e. $D:\Gamma^\infty(TP_0)\to \Gamma^\infty(TP_0)$ is linear and
\be\label{diffop2} D(f\mathcal{Z})=fD\mathcal{Z}+\mathfrak{X}(f)\mathcal{Z},\ee
where $f\in C^\infty(P_0)$. In order to obtain (\ref{diffop}) we note that
\be\label{fflow2} \frac{d}{dt}(\Gamma_t(\gamma_{-t}(\eta)))(w(\gamma_{-t}(\eta),\gamma_{-t}(\eta)^*))|_{t=0}= \ee $$=\frac{d}{dt}\Gamma_t(\gamma_{-t}(\eta))|_{t=0}(w(\eta,\eta^*))+\Gamma_0(\eta) \frac{d}{dt}w(\gamma_{-t}(\eta),\gamma_{-t}(\eta)^*)|_{t=0}=$$
$$=\frac{d}{dt}(\Gamma_{-t}(\eta)^{-1})|_{t=0}(w(\eta\eta^*))+\frac{d}{dt}w(\gamma_{-t}(\eta),\gamma_{-t}(\eta)^*)|_{t=0}=$$
$$=\Theta(\eta,\eta^*)(w(\eta,\eta^*))-\left\langle\frac{\partial w}{\partial\eta}(\eta,\eta^*),\vartheta(\eta,\eta^*)\right\rangle-\left\langle\frac{\partial w}{\partial\eta^*}(\eta,\eta^*),\vartheta(\eta,\eta^*)^*\right\rangle.$$
The bundle of derivations $\mathfrak{D}(TP_0)$ of the tangent bundle $TP_0$ is the algebroid of the trivial groupoid
$P_0\times G_0\times P_0\tto\Ll_{p_0}(\M)$ by virtue of the general theory developed for the finite dimensional case in Chapter 3 of \cite{mac}.

Eliminating $w:P_0\to \M p_0$ from (\ref{diffop}) we will obtain the following expression
\be\label{ddiffop3} D=-\mathfrak{X}+\Theta \ee
for the derivation $D\in \Gamma^\infty(\mathfrak{D}(TP_0))$. It is easy to see that the algebroid bracket of $D_1, D_2\in\Gamma^\infty(\mathfrak{D}(TP_0))$ is written in the form
\be\label{bracketkomut}[D_1,D_2]:=[\mathfrak{X}_1,\mathfrak{X}_2]+\mathfrak{X}_2(\Theta_1)-
\mathfrak{X}_1(\Theta_2)+[\Theta_1,\Theta_2]\ee

The anchor map $\mathbf{a}:\mathfrak{D}(TP_0)\to TP_0$ of the algebroid $\mathfrak{D}(TP_0)$ one defines by
\be\label{anchor} \mathbf{a}(D)=-\mathfrak{X}\ee
In particular when
\be\label{Gaamat} \Gamma_t(\eta)=T\gamma_t(\eta)\ee  the corresponding derivation  $D$ acts on $w$ as follows
\be\label{deriv}D(\eta)=\left(-\langle\frac{\partial w}{\partial\eta},\vartheta(\eta)\rangle-
\langle\frac{\partial w}{\partial\eta^*},\vartheta(\eta)^*
\rangle+\langle\frac{\partial \vartheta}{\partial\eta},w(\eta)\rangle+
\langle\frac{\partial \vartheta}{\partial\eta^*},w(\eta)^*\rangle,\eta\right).\ee
So, one has
\be D\mathcal{Z}=-[\mathfrak{X},\mathcal{Z}].\ee

The vector space consisting of derivations (\ref{deriv}) is closed with respect to the bracket (\ref{bracketkomut}). So, they form a Lie subalgebra of $\Gamma^\infty(\mathfrak{D}(TP_0))$ isomorphic to $\Gamma^\infty(TP_0)$. The isomorphism is given by the anchor map (\ref{anchor}).

According to (\ref{sigma3}) the one parameter subgroup of bisections $\{\sigma_t\}_{t\in \mathbb{R}} $ defines a one-parameter subgroup $\{\Phi_t\}_{t\in \mathbb{R}} $ of automorphisms $\Phi_t:TP_0\to TP_0$ of the tangent vector bundle. The vector field

\be\label{XD} \mathfrak{X}_D=\mathfrak{X}+\langle\frac{\partial }{\partial\vartheta},\Theta\vartheta\rangle+\langle\frac{\partial }{\partial\vartheta^*},(\Theta\vartheta)^*\rangle\ee
tangent to the flow $\{\Phi_t\}_{t\in \mathbb{R}} $ is a linear vector field on $TP_0$.

Since for $\mathfrak{X}_{D_1}, \mathfrak{X}_{D_2}\in \Gamma^\infty T(TP_0)$ one has $[\mathfrak{X}_{D_1}, \mathfrak{X}_{D_2}]=\mathfrak{X}_{[D_1,D_2]}$ we can identify the Lie algebra $\Gamma^\infty(\mathfrak{D}(TP_0))$ of derivations with the Lie algebra $\Gamma^\infty T(TP_0)$ of linear vector fields on the tangent vector bundle $TP_0$.

Now we are in a position  to investigate the algebroid $\A^{lin}_{p_0}(\M)$ of the groupoid $\G^{lin}\A_{p_0}(\M)\tto \Ll_{p_0}(\M)$ in terms of the algebroid $\mathfrak{D}(TP_0)$ of derivations of the tangent bundle $TP_0$. For this reason we note that any bisection of $\G^{lin}\A_{p_0}(\M)\tto \Ll_{p_0}(\M)$ is defined by the bisection (\ref{sigma}) as follows
\be\label{sigmaequi} [\sigma](l(\eta))[(\vartheta,\eta)]:=[(\Gamma(\eta)\vartheta,\gamma(\eta))],\ee
where the equivalence class $[\sigma(\eta)]:=[\sigma](l(\eta))$ is defined by the following equivalence relation:
\be\label{rel1}(\gamma(\eta),\Gamma(\eta),\eta)\sim (\gamma(\eta)h, R_h\circ\Gamma(\eta)\circ R_{g^{-1}},\eta g)\ee
for $h,g\in G_0$. From $l(\eta g)=l(\eta)$ taking into account  that (\ref{sigmaequi}) is independent of the choice of the classes representant we obtain the following condition on the bisection $\sigma$:
\be\label{sigmaequi2} \sigma(\eta g)=(\gamma(\eta g), \Gamma(\eta g), \eta g)=\ee
$$=(\gamma(\eta)h, R_h,\gamma(\eta))(\gamma(\eta),\Gamma(\eta), \eta)(\eta,R_{g^{-1}},\eta g)= $$
$$=(\gamma(\eta)h, R_h\circ\Gamma(\eta)\circ R_{g^{-1}},\eta g)$$
for any $h=h(\eta,g)\in G_0$ dependent on $\eta\in P_0$ and $g\in G_0$. We will assume that
$h:P_0\times G_0 \to  G_0$ is a smooth map. The condition (\ref{sigmaequi2}) is equivalent to the following conditions on $\gamma$ and $\Gamma$:
\be\label{cond1} \gamma(\eta g)=\gamma(\eta) h(\eta,g),\ee
\be\label{cond2} \Gamma(\eta g)=R_{h(\eta,g)}\circ\Gamma(\eta)\circ R_{g^{-1}},\ee
where $g\in G_0$. These conditions are consistent with the product formula (\ref{sigma2}). So the couple of maps $(\gamma, \Gamma)$ satisfying (\ref{cond1}) and (\ref{cond2}) defines bisection $[\sigma]$ of the groupoid $\G^{lin}\A_{p_0}(\M)\tto \Ll_{p_0}(\M)$ and any bisection of this groupoid is given in this way. From the above we also conclude that the product of $[\sigma_1]$ and $[\sigma_2]$ is correctly defined by
\be [\sigma_1]\ast[\sigma_2]:=[\sigma_1\ast\sigma_2].\ee
It follows from (\ref{cond1}) that $h(\eta, g)$ possesses the property of a cocycle
\be\label{cocycle} h(\eta, g_1 g_2)=h(\eta, g_1)h(\eta g_1, g_2).\ee
For the one-parameter group of bisections  $[\sigma_s](l(\eta))$ defined by $\sigma_s(\eta)=(\gamma_s(\eta),\Gamma_s(\eta),\eta)$ with exception  of the properties (\ref{bis1}) and (\ref{bis2}) we have the property
\be h_{s+t}(\eta,g)=h_s(\gamma_t(\eta), h_t(\eta,g)).\ee
Conditions (\ref{cond1}) and (\ref{cond2}) have their infinitesimal counterparts.
\be\label{condv} \vartheta(\eta g)=\vartheta(\eta) g+\eta H(\eta,g),\ee
\be\label{condteta}\Theta(\eta g)=R_g\circ \Theta(\eta)\circ R_{g^{-1}}+R{g^{-1}H(\eta,g)},\ee
where
\be H(\eta,g):=\frac{d}{dt}h_t(\eta,g)|_{t=0}.\ee

Summarizing the above facts we obtain:
\begin{prop} The algebroid $\A^{lin}_{p_0}(\M)$ could be considered as a subalgebroid of $\mathfrak{D}(TP_0)$ consisting of such derivations $D$, see (\ref{ddiffop3}), which satisfy conditions (\ref{condv}) and (\ref{condteta}). One concludes from this that these conditions are preserved by the bracket (\ref{bracketkomut}).

Further conditions
\be\label{thm} h_t(\eta,g)=g \quad { and}\quad \Gamma_t(\eta)=T\gamma_t(\eta)\ee
identify a subalgebroid of $\A^{lin}_{p_0}(\M)$ which is isomorphic to the algebroid $\A_{p_0}(\M)$ of the groupoid $\G_{p_0}(\M)\tto\Ll_{p_0}(\M)$.
\end{prop}

The results of this section one will find useful for investigation of various Poisson structures associated with groupoids $\G_{p_0}(\M)\tto\Ll_{p_0}(\M)$ and $\U_{p_0}(\M)\tto\Ll_{p_0}(\M)$. The paper concerning this subject is in preparation.

\end{document}